\documentclass[basic]{sn-jnl}
 \usepackage{times, mathptm}
 \usepackage{latexsym, amsmath, amssymb, bm, mathtools}
 \usepackage{hyperref, multirow}
 \usepackage[authoryear]{natbib}

\def\TS{\textstyle} \def\d{{\rm d}}
 \def\VEC#1{{\pmb{#1}}} \def\MAT#1{{\pmb{#1}}} \def\IE{{\it i.e.}} 
 \def\EG{{\it e.g.}}  
 \def\FIRST{{$1^{\textrm{st}}$}} \def\SECOND{{$2^{\textrm{nd}}$}} 
 \def\THIRD{{$3^{\textrm{rd}}$}} \def\NTH#1{{${#1}^{\textrm{th}}$}} 
 \def\DOTS{{...}}

 \def\AIJ{a_{ij}}
 \def\GAMMAC#1#2{\gamma_{{#1},c^{#2}}}
 \def\LAMBDAC#1#2{\lambda_{{#1}, c^{#2}}}
 \def\MUC#1#2{\mu_{{#1}, c^{#2}}}

 \def\STRUT{\vphantom{\vert_\vert^\vert}}
 \def\STRUTT{\vphantom{\Big\vert}}

\begin{document}

 \title{Embedded $(4, 5)$ pairs of explicit $7$-stage Runge--Kutta 
methods with FSAL property}

 \author*{\fnm{Misha} \sur{Stepanov}
 {\small ORCID: 0000-0002-7760-8239}
}\email{stepanov@math.arizona.edu}

 \affil{\orgdiv{Department of Mathematics \;and\; Program in Applied 
Mathematics}, \orgname{University of Arizona}, \orgaddress{
\city{Tucson}, \postcode{85721}, \state{AZ}, \country{USA}}}

\abstract{The general case of embedded $(4, 5)$ pairs of explicit 
$7$-stage Runge--Kutta methods with FSAL property ($a_{7 j} = b_{ j}$, 
$1 \le j \le 7$, $c_7 = 1$) is considered. Besides exceptional cases, 
the pairs form five $4$-dimensional families. The pairs within two 
(already known) families satisfy the simplifying assumption $\sum_j \AIJ 
c_{j} = c_i^2 / 2$, $i \ge 3$.}



\maketitle

\vspace{-24pt}

{\large\bf Declarations}

{\bf Conflicts of interest/Competing interests}: Not applicable.

\bigskip

Runge--Kutta methods (see, \EG, \citep[sec.~23 and ch.~3]{But16}, 
\citep[ch.~II]{HNW93}, \citep[ch.~4]{AsPe98}, \citep[ch.~3]{Ise08}) are 
widely and successfully used to solve Ordinary Differential Equations 
(ODEs) numerically for over a century \citep{BuWa96}. Consider a system 
$\d \VEC{x} / \d t = \VEC{f}(t, \VEC{x})$. To propagate by the step size 
$h$ and update the position, $\VEC{x}(t) \mapsto \VEC{\tilde{x}}(t + 
h)$, where $\VEC{\tilde{x}}(t + h)$ is a numerical approximation to the 
exact solution $\VEC{x}(t + h)$, an \mbox{$s$-stage} explicit 
Runge-Kutta method (which is determined by the coefficients $\AIJ$, 
weights $b_{ j}$, and nodes $c_{ i}$) would compute $\VEC{F}_1$, 
$\VEC{F}_2$, \DOTS, $\VEC{F}_{s}$, and then $\VEC{\tilde{x}}(t + 
h)$:\footnote{It is natural and will be assumed that $\sum_{j = 1}^{i - 
1} \AIJ = c_{ i}$. For $i = 1$ the sum is empty, so $c_1 = 0$ and 
$\VEC{F}_1 = \VEC{f} \bigl( t, \VEC{x}(t) \bigr)$.}
 \begin{align*} \VEC{F}_{i} = \VEC{f} \biggl( t + c_{ i} h, \VEC{x}(t) + 
h \sum_{j = 1}^{i - 1} \AIJ \VEC{F}_{j} \biggr), \qquad 
\VEC{\tilde{x}}(t + h) = \VEC{x}(t) + h \sum_{j = 1}^s b_{j} \VEC{F}_{j}
 \end{align*}

To obtain an accurate solution with less effort, various adaptive step 
size strategies were developed (see, \EG, \citep[sec.~33]{But16}, 
\citep[sec.~II.4]{HNW93}, \citep[sec.~4.5]{AsPe98}, 
\citep[ch.~6]{Ise08}). Typically the system of ODEs is solved in two 
different ways, and the step size is chosen so that the two solutions 
are sufficiently close. Embedded pairs of Runge--Kutta methods are 
computationally efficient, as the two methods within a pair have 
different weights, but share the nodes and the coefficients. The vectors 
$\VEC{F}_1$, $\VEC{F}_2$, \DOTS, $\VEC{F}_{s}$ are computed only once, 
and then are used in both methods.

The Butcher tableau \citep{But64} of an embedded $(4, 5)$ pair of 
explicit $7$-stage Runge--Kutta methods with so-called First Same As 
Last (FSAL) property \citep[p.~17]{Feh69}, \citep{DoPr78} looks like
 \begin{align*} \begin{array}{c|ccccccc}
     0\; \\
     c_2\; & \;a_{21}\; \\
     c_3\; & a_{31} & \;a_{32}\; \\
     c_4\; & a_{41} & a_{42} & \;a_{43}\; \\
     c_5\; & a_{51} & a_{62} & a_{53} & \;a_{54}\; \\
     c_6\; & a_{61} & a_{62} & a_{63} & a_{64} & \;a_{65}\; \\
     1\; & b_1 & b_2 & b_3 & b_4 & b_5 & b_6 \\
     \hline
         & b_1 & b_2 & b_3 & b_4 & b_5 & \;b_6\; \\
         & d_1 & d_2 & d_3 & d_4 & d_5 & d_6 & \;d_7\;
   \end{array} \end{align*} The vector $\VEC{b} = \bigl[ \, b_1 ~~ b_2 
~~ b_3 ~~ b_4 ~~ b_5 ~~ b_6 ~~ 0 \, \bigr]^{\textrm{T}}$ is the weights 
vector of the \NTH{5} order method, and $\VEC{d} = \bigl[ \, d_1 ~~ d_2 
~~ d_3 ~~ d_4 ~~ d_5 ~~ d_6 ~~ d_7 \, \bigr]^{\textrm{T}}$ is the 
difference between the \NTH{4} and the \NTH{5} order methods weights 
vectors.\footnote{Usually the \NTH{4} order method vector of 
weights $\VEC{b} + \VEC{d}$ is written in place of $\VEC{d}$.} The FSAL 
property means that the vector $\VEC{F}_1$ at the current step is equal 
to the already computed $\VEC{F}_7$ at the previous step. It implies 
$c_7 = 1$ and $a_{7  j} = b_{ j}$ for all $1 \le j 
\le 7$; \EG, $b_7 = 0$.

  \begin{table} \begin{tabular}{cccc}
    {\hskip0.2in{order}\hskip0.2in} & 
{\hskip0.45in{$\textrm{t}$}\hskip0.45in} & 
{\hskip0.45in{$\VEC{\Phi}(\textrm{t})$}\hskip0.45in} & 
{\hskip0.2in{$\textrm{t}!$}\hskip0.2in} \\
    \hline
    \FIRST{} & \begin{picture}(0,12)(0,2.5) \thicklines
  \put(0,6){\circle*{3}}
\end{picture} & $\VEC{1}$ & $1\STRUT$ \\
    \hline
    \SECOND{} & \begin{picture}(12,12)(0,3) \thicklines
  \put(0,6){\line(1,0){12}}
  \put(0,6){\circle*{3}}
  \put(12,6){\circle*{3}}
\end{picture} & $\VEC{c}$ & $2$ \\
    \hline
    \THIRD{} & \begin{picture}(12,12)(0,3) \thicklines
  \put(0,6){\line(4,-1){12}}
  \put(0,6){\line(4,1){12}}
  \put(0,6){\circle*{3}}
  \put(12,3){\circle*{3}}
  \put(12,9){\circle*{3}}
\end{picture} & $\VEC{c} * \VEC{c}$ & $3$ \\
    \raisebox{4pt}{$\downarrow$} &
\begin{picture}(24,12)(0,3) \thicklines
  \put(0,6){\line(1,0){24}}
  \put(0,6){\circle*{3}}
  \put(12,6){\circle*{3}}
  \put(24,6){\circle*{3}}
\end{picture} & $\VEC{c}'$ & $6\STRUT$ \\
    \hline
    \NTH{4} & \begin{picture}(12,12)(0,3) \thicklines
  \put(0,6){\line(1,0){12}}
  \put(0,6){\line(3,-1){12}}
  \put(0,6){\line(3,1){12}}
  \put(0,6){\circle*{3}}
  \put(12,2){\circle*{3}}
  \put(12,6){\circle*{3}}
  \put(12,10){\circle*{3}}
\end{picture} & $\VEC{c} * \VEC{c} * \VEC{c}$ & $4$ \\
  $\downarrow$ & \begin{picture}(24,12)(0,3) \thicklines
  \put(0,6){\line(4,-1){12}}
  \put(0,6){\line(4,1){12}}
  \put(12,3){\line(1,0){12}}
  \put(0,6){\circle*{3}}
  \put(12,3){\circle*{3}}
  \put(12,9){\circle*{3}}
  \put(24,3){\circle*{3}}
\end{picture} & $\VEC{c}' * \VEC{c}$ & $8$ \\
    & \begin{picture}(24,12)(0,3) \thicklines
  \put(0,6){\line(1,0){12}}
  \put(12,6){\line(4,-1){12}}
  \put(12,6){\line(4,1){12}}
  \put(0,6){\circle*{3}}
  \put(12,6){\circle*{3}}
  \put(24,3){\circle*{3}}
  \put(24,9){\circle*{3}}  
\end{picture} & $\MAT{A} (\VEC{c} * \VEC{c})$ & $12$ \\
    & \begin{picture}(36,12)(0,3) \thicklines
  \put(0,6){\line(1,0){36}}
  \put(0,6){\circle*{3}}
  \put(12,6){\circle*{3}}
  \put(24,6){\circle*{3}}
  \put(36,6){\circle*{3}}
\end{picture} & $\VEC{c}''$ & $24\STRUT$ \\
    \hline
    \NTH{5} & \begin{picture}(12,12)(0,3) \thicklines
  \put(0,6){\line(2,-1){12}}
  \put(0,6){\line(6,-1){12}}
  \put(0,6){\line(6,1){12}}
  \put(0,6){\line(2,1){12}}
  \put(0,6){\circle*{3}}
  \put(12,0){\circle*{3}}
  \put(12,4){\circle*{3}}
  \put(12,8){\circle*{3}}
  \put(12,12){\circle*{3}}
\end{picture} & $\VEC{c} * \VEC{c} * \VEC{c} * \VEC{c}$ & $5$ \\
    $\downarrow$ & \begin{picture}(24,12)(0,3) \thicklines
  \put(0,6){\line(1,0){12}}
  \put(0,6){\line(3,-1){12}}
  \put(0,6){\line(3,1){12}}
  \put(12,2){\line(1,0){12}}
  \put(0,6){\circle*{3}}
  \put(12,2){\circle*{3}}
  \put(12,6){\circle*{3}}
  \put(12,10){\circle*{3}}
  \put(24,2){\circle*{3}}
\end{picture} & $\VEC{c}' * \VEC{c} * \VEC{c}$ & $10$ \\
    & \begin{picture}(24,12)(0,3) \thicklines
  \put(0,7.5){\line(4,-1){12}}
  \put(0,7.5){\line(4,1){12}}
  \put(12,4.5){\line(4,-1){12}}
  \put(12,4.5){\line(4,1){12}}
  \put(0,7.5){\circle*{3}}
  \put(12,4.5){\circle*{3}}
  \put(12,10.5){\circle*{3}}
  \put(24,1.5){\circle*{3}}
  \put(24,7.5){\circle*{3}}
\end{picture} & $\bigl( \VEC{A} (\VEC{c} * \VEC{c}) \bigr) * \VEC{c}$ & 
$15$ \\
    & \begin{picture}(36,12)(0,3) \thicklines
  \put(0,6){\line(4,-1){12}}
  \put(0,6){\line(4,1){12}}
  \put(12,3){\line(1,0){24}}
  \put(0,6){\circle*{3}}
  \put(12,3){\circle*{3}}
  \put(12,9){\circle*{3}}
  \put(24,3){\circle*{3}}
  \put(36,3){\circle*{3}}
\end{picture} & $\VEC{c}'' * \VEC{c}$ & $30$ \\
    & \begin{picture}(24,12)(0,3) \thicklines
  \put(0,6){\line(4,-1){12}}
  \put(0,6){\line(4,1){12}}
  \put(12,3){\line(1,0){12}}
  \put(12,9){\line(1,0){12}}
  \put(0,6){\circle*{3}}
  \put(12,3){\circle*{3}}
  \put(12,9){\circle*{3}}
  \put(24,3){\circle*{3}}
  \put(24,9){\circle*{3}}
\end{picture} & $\VEC{c}' * \VEC{c}'$ & $20$ \\
    & \begin{picture}(24,12)(0,3) \thicklines
  \put(0,6){\line(1,0){12}}
  \put(12,6){\line(1,0){12}}
  \put(12,6){\line(3,-1){12}}
  \put(12,6){\line(3,1){12}}
  \put(0,6){\circle*{3}}
  \put(12,6){\circle*{3}}
  \put(24,2){\circle*{3}}
  \put(24,6){\circle*{3}}
  \put(24,10){\circle*{3}}
\end{picture} & $\MAT{A} (\VEC{c} * \VEC{c} * \VEC{c})$ & $20$ \\
    & \begin{picture}(36,12)(0,3) \thicklines
  \put(0,6){\line(1,0){12}}
  \put(12,6){\line(4,-1){12}}
  \put(12,6){\line(4,1){12}}
  \put(24,3){\line(1,0){12}}
  \put(0,6){\circle*{3}}
  \put(12,6){\circle*{3}}
  \put(24,3){\circle*{3}}
  \put(24,9){\circle*{3}}
  \put(36,3){\circle*{3}}
\end{picture} & $\MAT{A} (\VEC{c}' * \VEC{c})$ & $40$ \\
    & \begin{picture}(36,12)(0,3) \thicklines
  \put(0,6){\line(1,0){24}}
  \put(24,6){\line(4,-1){12}}
  \put(24,6){\line(4,1){12}}
  \put(0,6){\circle*{3}}
  \put(12,6){\circle*{3}}
  \put(24,6){\circle*{3}}
  \put(36,3){\circle*{3}}
  \put(36,9){\circle*{3}}
\end{picture} & $\MAT{A}^2 (\VEC{c} * \VEC{c})$ & $60$ \\
    & \begin{picture}(48,12)(0,3) \thicklines
  \put(0,6){\line(1,0){48}}
  \put(0,6){\circle*{3}}
  \put(12,6){\circle*{3}}
  \put(24,6){\circle*{3}}
  \put(36,6){\circle*{3}}
  \put(48,6){\circle*{3}}
\end{picture} & $\VEC{c}'''$ & $120\STRUT$
  \end{tabular}
  \caption{Order conditions $\VEC{b}^{\rm T} \VEC{\Phi}(\textrm{t}) = 1 / 
\textrm{t}!$ for rooted trees $\textrm{t}$ with $1$, $2$, $3$, $4$, and 
$5$ vertices. The ``$*$'' sign denotes component-wise multiplication of 
vectors, \IE, $\bigl( \VEC{x} * \VEC{y} \bigr)_i = x_i y_i$.} 
\label{order_conditions} \end{table}

Let $\VEC{c} = \bigl[ \, 0 ~~ c_2 ~~ c_3 ~~ c_4 ~~ c_5 ~~ c_6 ~~ 1 \, 
\bigr]^{\textrm{T}}$ and $\MAT{A} = \bigl[ \AIJ \bigr]$ be the $7 \times 
7$ matrix with $\AIJ$ as its matrix element in the \NTH{i} row and 
\NTH{j} column (visibly, $\AIJ = 0$ if $i \le j$). Let $\VEC{1}$ be the 
vector with all components being equal to $1$. The condition $\sum_j 
\AIJ = c_{ i}$ or $\MAT{A} \VEC{1} = \VEC{c}$ is assumed. Let $\VEC{c}' 
= \MAT{A} \VEC{c}$, $\VEC{c}'' = \MAT{A} \VEC{c}'$, and $\VEC{c}''' = 
\MAT{A} \VEC{c}'' = \MAT{A}^4 \VEC{1}$. A Runge--Kutta method of order 
$p$ should satisfy the conditions $\VEC{b}^{\rm T} 
\VEC{\Phi}(\textrm{t}) = 1 / \textrm{t}!$ for all rooted trees 
$\textrm{t}$ with up to $p$ vertices \citep[p.~175]{But16}, 
\citep[p.~177]{But21}, \citep[p.~153]{HNW93}. For $p = 5$ these 
conditions are listed in Table~\ref{order_conditions}, see also 
\citep[p.~172]{But16}, \citep[p.~126]{But21}, \citep[p.~148]{HNW93}, 
\citep[tab.~1]{DoPr80}. For a $(4, 5)$ pair the conditions $\VEC{b}^{\rm 
T} \VEC{\Phi}(\textrm{t}) = 1 / \textrm{t}!$ and $\VEC{d}^{\rm T} 
\VEC{\Phi}(\textrm{t}) = 0$ are satisfied for all trees $\textrm{t}$ 
with up to $5$ and $4$ vertices, respectively.

The process of Runge--Kutta methods construction is streamlined by using 
so-called simplifying assumptions (see, \EG, \citep[sec.~321]{But16}, 
\citep[pp.~136 and 175]{HNW93}). The one that is important for the 
subject discussed here is $c'_i = \sum_j a_{ij} c_j = c_i^2 / 2$ for any 
$i \ne 2$. For a method of order at least $3$ this would imply $b_2 = 
0$, as $\VEC{b}^{\textrm{T}} \VEC{c}' = \frac{1}{6} = \frac{1}{2} 
\VEC{b}^{\textrm{T}} (\VEC{c} * \VEC{c})$.

There is no $5$-stage explicit Runge--Kutta \NTH{5} order method 
\citep{But64}. The general case of the $6$-stage, \NTH{5} order method 
was considered in \citep{Cas66}, \citep{Cas69}, where the set of order 
conditions, by exclusion of variables, was drastically reduced, and 
methods with $b_2 \ne 0$ were built. A two-dimensional family of 
embedded $(4, 5)$ pairs of $6$-stage Runge--Kutta methods (with $3 c_2 = 
2 c_3$, $c_5 = 1$, and $b_2 = d_2 = d_7 = 0$) was constructed in 
\citep{Feh69}. A method suggested in \citep{CaKa90} belongs to this 
family. In \citep{DoPr80} a three-dimensional family of $7$-stage pairs 
with FSAL property was presented (with $3 c_2 = 2 c_3$, $c_6 = 1$, and 
$b_2 = d_2 = 0$). Both families were extended to four-dimensional ones 
in \citep{PaPa96}. In \citep{Tsi11} FSAL pairs not satisfying the 
simplifying assumption were considered, seemingly with the aim of 
extending the set of pairs satisfying the order conditions and thus 
potentially finding a more efficient and practical pair. With the 
conditions solved in part analytically and in part numerically, the 
\citep[tab.~1]{Tsi11} pair was suggested.

Increasing the number of stages (and thus the amount of computation per 
step) provides additional flexibility in choosing $\MAT{A}$, $\VEC{b}$, 
$\VEC{c}$, and $\VEC{d}$, which may be exploited to construct viable 
pairs that produce an accurate solution in fewer steps. In 
\citep[sec.~3.1]{ShSm93} and \citep{BoSh96} non-FSAL embedded $(4, 5)$ 
pairs of $7$-stage Runge--Kutta methods were suggested.

In this paper embedded $(4, 5)$ pairs of $7$-stage Runge--Kutta methods 
with FSAL property (this includes non-FSAL pairs of $6$-stage methods) 
are considered, with the aim of complete classification of at least 
general, non-exceptional, cases. After rewriting the order conditions in 
terms of $\VEC{c}$, $\VEC{c}'$, $\VEC{c}''$, $\VEC{c}'''$, and $a_{65}$, 
$b_6$, $d_5$, $d_6$, $d_7$ (Section~\ref{rewrite}), a pair is expressed 
through $6$ variables: $c_2$, $c_3$, $c_4$, $c_5$, $c_6$, and $c'_3$ 
(Section~\ref{express}). Lastly, pairs are classified into five 
$4$-dimensional families (Section~\ref{five_families}). The topics of 
choosing the magnitude of the vector $\VEC{d}$ and of continuous 
formulas or interpolants (see, \EG, \citep[sec.~II.6]{HNW93}) are not 
considered.

\section{Rewriting some of the order conditions in a compact form} 
\label{rewrite}

\hskip\parindent{}It is convenient to express the $30 = 5_{\VEC{c} = 
\MAT{A} \VEC{1}} + 17_{\textrm{\NTH{5}~order,~}\VEC{b}} + 
8_{\textrm{\NTH{4}~order,~}\VEC{d}}$ conditions on an embedded pair not 
in terms of $\MAT{A}$, $\VEC{b}$, $\VEC{c}$, and $\VEC{d}$ ($33 = 
21_{\MAT{A}, \VEC{b}} + 5_{\VEC{c}} + 7_{\VEC{d}}$ degrees of freedom), 
but in terms of $\VEC{c}$, $\VEC{c}'$, $\VEC{c}''$, $\VEC{c}'''$, 
$a_{65}$, $b_6$, $d_5$, $d_6$, and $d_7$ ($19 = 5_{\VEC{c}} + 
4_{\VEC{c}'} + 3_{\VEC{c}''} + 2_{\VEC{c}'''} + 5_{a_{65},  
b_6,  d_{5, 6, 7}}$ degrees of freedom). After this (rather 
mechanical) change of variables the relations $\VEC{c} = \MAT{A} 
\VEC{1}$, $\VEC{c}' = \MAT{A} \VEC{c}$, $\VEC{c}'' = \MAT{A} \VEC{c}'$, 
$\VEC{c}''' = \MAT{A} \VEC{c}''$ and the order conditions 
$\VEC{b}^{\textrm{T}} \bigl[ \, \VEC{1} ~~ \VEC{c} ~~ \VEC{c}' ~~ 
\VEC{c}'' ~~ \VEC{c}''' \, \bigr] = \bigl[ \, 1 ~~ \frac{1}{2} ~~ 
\frac{1}{6} ~~ \frac{1}{24} ~~ \frac{1}{120} \, \bigr]$, 
$\VEC{d}^{\textrm{T}} \bigl[ \, \VEC{1} ~~ \VEC{c} ~~ \VEC{c}' ~~ 
\VEC{c}'' \, \bigr] = \bigl[ \, 0 ~~ 0 ~~ 0 ~~ 0 \, \bigr]$ will be 
satisfied by construction.\footnote{The FSAL property $a_{7 
 j} = b_{ j}$, $1 \le j \le 7$ and the order 
conditions $\VEC{b}^{\textrm{T}} \VEC{c} = 1 / 2$, $\VEC{b}^{\textrm{T}} 
\VEC{c}' = 1 / 6$, $\VEC{b}^{\textrm{T}} \VEC{c}'' = 1 / 24$ result in 
$c'_7 = 1 / 2$, $c''_7 = 1 / 6$, and $c'''_7 = 1 / 24$.} There still 
going to be $16 = (17 - 5)_{\textrm{\NTH{5}~order,~}\VEC{b}} + (8 - 
4)_{\textrm{\NTH{4}~order,~}\VEC{d}}$ (redundant) order conditions left.

The condition $\VEC{c} = \MAT{A} \VEC{1}$ and the order conditions 
$\VEC{b}^{\textrm{T}} \VEC{1} = 1$, $\VEC{d}^{\textrm{T}} \VEC{1} = 0$ 
imply
 \begin{align*}
   a_{21} &= c_2 \\
   a_{31} &= c_3 - a_{32} \\
   a_{41} &= c_4 - a_{42} - a_{43} \\
   a_{51} &= c_5 - a_{52} - a_{53} - a_{54} \\ 
   a_{61} &= c_6 - a_{62} - a_{63} - a_{64} - a_{65} \\
   b_1 &= 1 - b_2 - b_3 - b_4 - b_4 - b_6 \\
   d_1 &= -d_2 - d_3 - d_4 - d_4 - d_6 - d_7 \qquad\quad\;\,\,
 \end{align*} Five stages are not enough to satisfy all the required 
order conditions \citep{But64}, thus $c_2 \ne 0$ (otherwise the \FIRST{} 
and \SECOND{} stages are redundant) and $b_6 \ne 0$. The relation 
$\VEC{c}' = \MAT{A} \VEC{c}$ and the order conditions 
$\VEC{b}^{\textrm{T}} \VEC{c} = \frac{1}{2}$, $\VEC{d}^{\textrm{T}} 
\VEC{c} = 0$ imply
 \begin{align*}
   a_{32} &= c'_3 / c_2 \\
   a_{42} &= (c'_4 - a_{43} c_3) / c_2 \\
   a_{52} &= (c'_5 - a_{53} c_3 - a_{54} c_4) / c_2 \\
   a_{62} &= (c'_6 - a_{63} c_3 - a_{64} c_4 - a_{65} c_5) / c_2 \\
      b_2 &= (1 / 2 - b_3 c_3 - b_4 c_4 - b_5 c_5 - b_6 c_6) / c_2 \\
      d_2 &= (-d_3 c_3 - d_4 c_4 - d_5 c_5 - d_6 c_6 - d_7) / c_2 
 \end{align*} In what follows it is going to be assumed that the matrix 
elements of $\MAT{A}$ right below the diagonal are non-zero: $a_{32} \ne 
0$, $a_{43} \ne 0$, $a_{54} \ne 0$, and $a_{65} \ne 0$.\footnote{The 
full analysis of $a_{65} a_{54} a_{43} a_{32} = 0$ case is tedious and 
is not expected to result in an embedded pair of practical interest. For 
instance, if $a_{32} = 0$, then $c_3 = 3 c_2 / (8 c_2 - 3)$ and $c_4 = 
0$.} This is equivalent to $c'_3 \ne 0$, $c''_4 \ne 0$, $c'''_5 \ne 0$, 
and $a_{65} \ne 0$. The relations $\VEC{c}'' = \MAT{A} \VEC{c}'$, 
$\VEC{c}''' = \MAT{A} \VEC{c}''$ and the order conditions 
$\VEC{b}^{\textrm{T}} \bigl[ \, \VEC{c}' ~~ \VEC{c}'' ~~ \VEC{c}''' \, 
\bigr] = \bigl[ \, \frac{1}{6} ~~ \frac{1}{24} ~~ \frac{1}{120} \, 
\bigr]$, $\VEC{d}^{\textrm{T}} \bigl[ \, \VEC{c}' ~~ \VEC{c}'' \, \bigr] 
= \bigl[ \, 0 ~~ 0 \, \bigr]$ imply
 \begin{align*}
   a_{43} &= c''_4  / c'_3 \\
   a_{53} &= (c''_5 - a_{54} c'_4) / c'_3 \\
   a_{63} &= (c''_6 - a_{64} c'_4 - a_{65} c'_5) / c'_3 \\
      b_3 &= (1 / 6 - b_4 c'_4 - b_5 c'_5 - b_6 c'_6) / c'_3 \\
      d_3 &= (-d_4 c'_4 - d_5 c'_5 - d_6 c'_6 - d_7 / 2) / c'_3 \\
   a_{54} &= c'''_5 / c''_4 \\
   a_{64} &= (c'''_6 - a_{65} c''_5) / c''_4 \\
      b_4 &= (1 / 24 - b_5 c''_5 - b_6 c''_6) / c''_4 \\
      d_4 &= (-d_5 c''_5 - d_6 c''_6 - d_7 / 6) / c''_4 \\
      b_5 &= (1 / 120 - b_6 c'''_6) / c'''_5
 \end{align*} Now $a_{ i  j}$, $b_{ j}$, 
$d_{ j}$, where $2 \le i \le 7$, $1 \le j \le 4$, and $b_5$ 
are expressed through $\VEC{c}$, $\VEC{c}'$, $\VEC{c}''$, $\VEC{c}'''$, 
$a_{65}$, $b_6$, $d_5$, $d_6$, and $d_7$. The variables $b_5$ and 
$c'''_6$ are interchangeable:
 \begin{gather} b_5 = (1 / 120 - b_6 c'''_6) / c'''_5 \quad 
\longleftrightarrow \quad c'''_6 = (1 / 120 - b_5 c'''_5) / b_6 
\label{b5_cppp6} \\ a_{64} = (c'''_6 - a_{65} c''_5) / c''_4 = (1 / 120 
- b_5 c'''_5) / b_6 c''_4 - a_{65} c''_5 / c''_4 \nonumber \end{gather}

The following notation will be useful, where $4 \le m \le 7$ and $1 \le 
n \le 3$:\footnote{Further derivation was done in interaction with 
computer algebra system Wolfram Mathematica~8.0, mainly using commands 
\scalebox{0.725}[0.825]{\scriptsize\tt\bfseries Solve} to symbolically 
solve linear equations, \scalebox{0.725}[0.825]{\scriptsize\tt\bfseries 
Simplify}, and (in Section~\ref{five_families}) 
\scalebox{0.725}[0.825]{\scriptsize\tt\bfseries Factor}.}
 \begin{align*}
   \gamma_{m,  c^{n + 1}} &= c'_3 c_m (c_m^n - c_2^n) - c'_m 
c_3 (c_3^n - c_2^n) \\
   \gamma_{m,  c'  c^n} &= c'_m (c_m^n - c_3^n) \\
   \gamma_{m,  c'^2} &= c'_m (c'_m - c'_3) \\
   \gamma_{m,  c''  c} &= c''_m (c_m - c_4) \\
   \lambda_{m, *} &= c''_4 \gamma_{m, *} - c''_m \gamma_{4, *}, \qquad 
\smash{* = c^{n + 1},  c'  c^n,  c'^2, 
 c''  c} \\
   \gamma_{ A c^2} &= c''_4 c_3 (c_3 - c_2) \\
   \mu_{m,  c^{n + 1}} &= \gamma_{m, c^{n + 1}} + 4 c''_m 
\bigl( c_3 (c_3^n - c_2^n) + 3 c'_3 (c_2^n - \smash{\TS\frac{2}{n + 2}}) 
\bigr) \\
   \mu_{m,  c'  c^n} &= \gamma_{m, c^n c'} + 4 c''_m 
(c_3^n - \smash{\TS\frac{3}{n + 3}}) \\
   \mu_{m,  c'^2} &= \gamma_{m, c'^2} + 4 c''_m (c'_3 - 
\smash{\TS\frac3{10}}) \\
   \mu_{m,  c''  c} &= c''_m (c_m - 
\smash{\TS\frac45}) \\
   \eta_{c^{n + 1}} &= c_3 (c_3^n - c_2^n) + 4 c'_3 (c_2^n - 
{\TS\frac3{4 n + 2}}) \\
   \eta_{c'  c} &= c_3 - \smash{\TS\frac35} \\
   \eta_{A c^2} &= c'_3 (c_2 - \smash{\TS\frac25})
 \end{align*}

The remaining $12$ order conditions for the \NTH{5}{} order method, with 
the exception of $\VEC{b}^{\textrm{T}} \bigl( (\MAT{A} (\VEC{c} * 
\VEC{c})) * \VEC{c} \bigr) = \frac{1}{15}$, could be written as 
\begin{align} {\rm rank} \left[ \begin{array}{cc} b_6 & 1 / 120 \\
  \LAMBDAC{5}{2} + 5 c'''_5 \MUC{4}{2} & c'''_6 \LAMBDAC{5}{2} - c'''_5 
\LAMBDAC{6}{2} \\
  \LAMBDAC{5}{3} + 5 c'''_5 \MUC{4}{3} & c'''_6 \LAMBDAC{5}{3} - c'''_5 
\LAMBDAC{6}{3} \\
  \LAMBDAC{5}{4} + 5 c'''_5 \MUC{4}{4} & c'''_6 \LAMBDAC{5}{4} - c'''_5 
\LAMBDAC{6}{4} \\
  \lambda_{5,  c'  c} + 5 c'''_5 \mu_{4,  c' 
 c} & c'''_6 \lambda_{5,  c'  c} - c'''_5 
\lambda_{6,  c'  c} \\
  \lambda_{5,  c'  c^2} + 5 c'''_5 \mu_{4,  
c'  c^2} & c'''_6 \lambda_{5,  c'  c^2} - 
c'''_5 \lambda_{6,  c'  c^2} \\
  \lambda_{5, c'^2} + 5 c'''_5 \mu_{4, c'^2} & c'''_6 \lambda_{5, c'^2} 
- c'''_5 \lambda_{6, c'^2} \\
  \lambda_{5, c'' c} + 5 c'''_5 \mu_{4, c'' c} & c'''_6 \lambda_{5, c'' 
c} - c'''_5 \lambda_{6, c'' c} \\
  \GAMMAC{4}{2} + 5 c''_4 \eta_{c^2} & -a_{65} \LAMBDAC{5}{2} \\ 
  \GAMMAC{4}{3} + 5 c''_4 \eta_{c^3} & -a_{65} \LAMBDAC{5}{3} \\ 
  \gamma_{4,  c'  c} + 5 c''_4 \eta_{c'  c} 
& -a_{65} \lambda_{5,  c'  c} \\
  \gamma_{A c^2} + 5 c''_4 \eta_{A c^2} & -a_{65} c'''_5 \GAMMAC{4}{2} 
  \end{array} \right] = 1 \label{b6_12x2} \end{align} Currently the 
whole vector $\VEC{b}$ is expressed through $\VEC{c}$, $\VEC{c}'$, 
$\VEC{c}''$, $\VEC{c}'''$, $a_{65}$, and $b_6$. Any but the \FIRST{} row 
in this $12 \times 2$ matrix gives the solution for $b_6$ in the 
corresponding order condition $\VEC{b}^{\textrm{T}} 
\VEC{\Phi}(\textrm{t}) = 1 / \textrm{t}!$. From the second to eighth row 
these conditions can be rewritten as\footnote{Also $b_4 \mu_{4, *} + b_5
 \mu_{5, *} + b_6  \mu_{6, *} = 0$, as $b_4 = (1 / 24 
- b_5 c''_5 - b_6 c''_6) / c''_4$.} \begin{align} \hskip-84pt 
\overbrace{\left[ \begin{array}{ccc} {\mu_{4, *}} & {\mu_{5, *}} 
& {\mu_{6, *}} \end{array} \right] \left[ \begin{array}{ccc} 1 & 
-c''_5 & -c''_6 \\ 0 & \phantom{-}c''_4 & ~0 \\ 0 & ~0 & 
\phantom{-}c''_4 \end{array} \right]}^{\small \left[ \begin{array}{ccc} 
{\mu_{4, *}} & {\lambda_{5, *}} & {\lambda_{6, *}} 
\end{array} \right]} \overbrace{\left[ \begin{array}{cc} {5 c'''_5} 
& 0 \\ 1 & {\phantom{-}c'''_6} \\ 0 & {-c'''_5} \end{array} 
\right] \left[ \begin{array}{c} {1 / 120} \\ {-b_6} \end{array} 
\right]}^{\small c'''_5  \left[ \begin{array}{ccc} {1 / 24} 
& {b_5} & {b_6} \end{array} 
\right]^{\textrm{\scriptsize T}}} = \left[ \begin{array}{c} 
{0} \end{array} \right] \hskip-257.5pt \underbrace{\phantom{\left[ 
\begin{array}{ccc} {\mu_{4, *}} & {\mu_{5, *}} & {\mu_{6, 
*}} \end{array} \right] \left[ \begin{array}{ccc} 1 & -c''_5 & -c''_6 
\\ 0 & \phantom{-}c''_4 & ~0 \\ 0 & ~0 & \phantom{-}c''_4 \end{array} 
\right] \left[ \begin{array}{cc} {5 c'''_5} & 0 \\ 1 & 
{\phantom{-}c'''_6} \\ 0 & {-c'''_5} \end{array} 
\right]}}_{\small \left[ \begin{array}{cc} {\lambda_{5, *} + 5 c'''_5 
\mu_{4, *}} & {c'''_6 \lambda_{5, *} - c'''_5 \lambda_{6, *}} 
\end{array} \right]} \label{mu456} \end{align}

For the \NTH{4} and \NTH{5} order methods in the pair to produce 
distinct solutions, the vector $\VEC{d}$ is non-zero. The following four 
combinations should be equal to zero:
 \begin{align*}
  c'_3 c''_4 \, \VEC{d}^{\rm T} (\VEC{c} * \VEC{c}) &= d_5 \lambda_{5, 
 c^2} + d_6 \LAMBDAC{6}{2} + d_7 \lambda_{7, 
 c^2} \hskip1.3in \scalebox{0.7}{\begin{picture}(12,12)(0,3) 
\thicklines
  \put(0,6){\line(4,-1){12}}
  \put(0,6){\line(4,1){12}}
  \put(0,6){\circle*{3}}
  \put(12,3){\circle*{3}}
  \put(12,9){\circle*{3}}
\end{picture}}
 \\
  c'_3 c''_4 \, \VEC{d}^{\rm T} (\VEC{c} * \VEC{c} * \VEC{c}) &= d_5 
\LAMBDAC{5}{3} + d_6 \LAMBDAC{6}{3} + d_7 
\LAMBDAC{7}{3} \hskip1.3in 
\scalebox{0.7}{\begin{picture}(12,12)(0,3) \thicklines
  \put(0,6){\line(1,0){12}}
  \put(0,6){\line(3,-1){12}}
  \put(0,6){\line(3,1){12}}
  \put(0,6){\circle*{3}}
  \put(12,2){\circle*{3}}
  \put(12,6){\circle*{3}}
  \put(12,10){\circle*{3}}
\end{picture}}
 \\
  c''_4 \, \VEC{d}^{\rm T} (\VEC{c}' * \VEC{c}) &= d_5 \lambda_{5, 
 c'  c} + d_6 \lambda_{6, c' c} 
+ d_7 \lambda_{7,  c'  c} \hskip1.135in 
\scalebox{0.7}{\begin{picture}(24,12)(0,3) \thicklines
  \put(0,6){\line(4,-1){12}}
  \put(0,6){\line(4,1){12}}
  \put(12,3){\line(1,0){12}}
  \put(0,6){\circle*{3}}
  \put(12,3){\circle*{3}}
  \put(12,9){\circle*{3}}
  \put(24,3){\circle*{3}}
\end{picture}}
 \\
  c'_3 c''_4 \, \VEC{d}^{\rm T} \MAT{A} (\VEC{c} * \VEC{c}) &= d_5 
c'''_5 \GAMMAC{4}{2} + d_6 (c'''_6 \GAMMAC{4}{2} 
+ a_{65} \LAMBDAC{5}{2}) \hskip0.8in 
\raisebox{1pt}{\scalebox{0.7}{\begin{picture}(24,12)(0,3) \thicklines
  \put(0,6){\line(1,0){12}}
  \put(12,6){\line(4,-1){12}}
  \put(12,6){\line(4,1){12}}
  \put(0,6){\circle*{3}}
  \put(12,6){\circle*{3}}
  \put(24,3){\circle*{3}}
  \put(24,9){\circle*{3}}
\end{picture}}} \\
  &+ d_7 \bigl( \LAMBDAC{5}{2} + 5 c'''_5 \gamma_{4, 
 c^2} - 120 b_6 (c'''_6 \LAMBDAC{5}{2} - c'''_5 
\LAMBDAC{6}{2}) \bigr) / 120 c'''_5
 \end{align*} The condition $\VEC{b}^{\rm T} (\VEC{c} * \VEC{c}) = 
\frac{1}{3}$ implies $120 b_6 (c'''_6 \LAMBDAC{5}{2} - 
c'''_5 \LAMBDAC{6}{2}) = \LAMBDAC{5}{2} + 5 
c'''_5 \MUC{4}{2}$, which simplifies the coefficient at 
$d_7$ in $\VEC{d}^{\rm T} \MAT{A} (\VEC{c} * \VEC{c})$. As the 
conditions on the vector $\VEC{d}$ are linear and homogeneous, it can be 
rescaled by any non-zero factor. Such a rescaling just recalibrates the 
measure of closeness between the two solutions in the adaptive step size 
scheme.

Here are the conditions on the vector $\VEC{d}$ combined with 
eq.~(\ref{mu456}) for $* = c^2$, $c^3$, and $c' c$, and also with the 
last four rows of the matrix in 
eq.~(\ref{b6_12x2}):\footnote{In the case of an embedded pair 
of $6$-stage Runge--Kutta methods, \IE, $d_7 = 0$, the rank of the 
matrix $\MAT{M}$ without the \NTH{5} row should be equal to $1$.}
 \begin{gather} \left[ \begin{array}{ccccc}
      {1/ 120} & 0            & {b_6 a_{65}} & 0   & 0 \\
      0            & {1 / 24} & b_5              & b_6 & 0 \\
      0            & 0            & d_5              & d_6 & d_7
    \end{array} \right] \MAT{M} =
    \left[ \begin{array}{cccc}
      0 & 0 & 0 & 0 \\
      0 & 0 & 0 & 0 \\
      0 & 0 & 0 & 0 \end{array} \right] \label{ABDM} \\
  \MAT{M} = \left[ \begin{array}{cccc}
   {\GAMMAC{4}{2} + 5 c''_4 \eta_{c^2}} & {\gamma_{4, 
 c^3} + 5 c''_4 \eta_{c^3}} & {\gamma_{4,  c' 
 c} + 5 c''_4 \eta_{c'  c}} & {\gamma_{A c^2} + 
5 c''_4 \eta_{A c^2}} \\
   \MUC{4}{2} & \MUC{4}{3} & \mu_{4,  
c'  c} & c''_4 \eta_{c^2} \\
   \LAMBDAC{5}{2} & \LAMBDAC{5}{3} & \lambda_{5, 
 c'  c} & c'''_5 \GAMMAC{4}{2} \\
   \LAMBDAC{6}{2} & \LAMBDAC{6}{3} & \lambda_{6,  
c'  c} & c'''_6 \GAMMAC{4}{2} + a_{65} 
\LAMBDAC{5}{2} \\
   \LAMBDAC{7}{2} & \LAMBDAC{7}{3} & \lambda_{7,  
c'  c} & (\GAMMAC{4}{2} - \MUC{4}{2}) 
/ 24
 \end{array} \right] \nonumber \end{gather} As $c''_4 \eta_{c^2} / 24 + 
b_5 c'''_5 \GAMMAC{4}{2} + b_6 (c'''_6 \gamma_{4,  
c^2} + a_{65} \LAMBDAC{5}{2}) = (\GAMMAC{4}{2} + 
5 c''_4 \eta_{c^2}) / 120 + b_6 a_{65} \LAMBDAC{5}{2} = 0$, 
the matrix element in the \SECOND{} row and the \NTH{4} column of the 
product in eq.~(\ref{ABDM}) is equal to zero. All the columns of 
$\MAT{M} = \bigl[ m_{ i  j} \bigr]$ are orthogonal to 
any row of the $3 \times 5$ matrix in eq.~(\ref{ABDM}) whose rank is 
$3$, thus $\textrm{rank\,}\MAT{M} \le 2$. The \FIRST{} and \THIRD{} rows 
of $\MAT{M}$ are proportional to each other.

The order conditions that are not taken into account in eq.~(\ref{ABDM}) 
are $\VEC{b}^{\textrm{T}} \bigl[ \, (\VEC{c} * \VEC{c} * \VEC{c} * 
\VEC{c}) ~~~ (\VEC{c}' * \VEC{c} * \VEC{c}) ~~~ \bigl( \bigl( \MAT{A} 
(\VEC{c} * \VEC{c}) \bigr) * \VEC{c} \bigr) ~~~ (\VEC{c}'' * \VEC{c}) 
~~~ (\VEC{c}' * \VEC{c}') \, \bigr] = \bigl[ \, \frac{1}{5} ~~ 
\frac{1}{10} ~~ \frac{1}{15} ~~ \frac{1}{30} ~~ \frac{1}{20} \, \bigr]$.

\section{Expressing a pair through $c_2$, $c_3$, $c_4$, $c_5$, $c_6$, 
and $c'_3$} \label{express}

The first two rows of eq.~(\ref{ABDM}) are used to express $a_{65}$, 
$b_5$, $b_6$, and $c'''_5$ (and thus also $c'''_6$, see 
eq.~(\ref{b5_cppp6})) through $\VEC{c}$, $\VEC{c}'$, and $\VEC{c}''$:
 \begin{align*}
   b_5 = -\frac{1}{24} \frac {\MUC{4}{2} \lambda_{6, 
 c'  c} - \mu_{4, c' c} \LAMBDAC{6}{2}} 
{\LAMBDAC{5}{2} \lambda_{6,  c'  c} - \lambda_{5, 
 c'  c} \LAMBDAC{6}{2}} , & \qquad
   c'''_5 = \frac {m_{31} m_{14}} {\GAMMAC{4}{2} m_{11}} = \frac 
{\LAMBDAC{5}{2} (\gamma_{A c^2} + 5 c''_4 \eta_{A c^2})} {\gamma_{4, 
c^2} (\GAMMAC{4}{2} + 5 c''_4 \eta_{c^2})} \\
   b_6 = \frac{1}{24} \frac {\MUC{4}{2} \lambda_{5, 
 c'  c} - \mu_{4,  c'  c} 
\LAMBDAC{5}{2}} {\LAMBDAC{5}{2} \lambda_{6,  c'  c} - 
\lambda_{5,  c'  c} \LAMBDAC{6}{2}} , & \qquad
   a_{65} = - \frac{m_{11}}{120 b_6 m_{31}} = - \frac {\GAMMAC{4}{2} + 
5 c''_4 \eta_{c^2}} {120 b_6 \lambda_{5, c^2}}
 \end{align*} The element $m_{44}$ of the matrix $\MAT{M}$ is equal to 
$m_{44} = c'''_6 \GAMMAC{4}{2} + a_{65} \LAMBDAC{5}{2} = 
\GAMMAC{4}{2} \bigl( \frac{1}{120} - b_5 c'''_5 \bigr) / b_6 
- (\GAMMAC{4}{2} + 5 c''_4 \eta_{c^2}) / 120 b_6 = -(m_{24} 
/ 24 + b_5 m_{34}) / b_6$, which is compatible with the \SECOND{} row of 
eq.~(\ref{ABDM}).

By performing the following elementary row and column operations (the 
order is important) the matrix $\MAT{M}$ is brought to a simpler form:
 \begin{align*}
   \MAT{M}_{*2} &\leftarrow \MAT{M}_{*2} - (c_2 + c_3) \MAT{M}_{*1} \\
   \MAT{M}_{(m - 2)*} &\leftarrow \bigl( \MAT{M}_{(m - 2)*} + c''_m 
 \MAT{M}_{1*} \bigr) / c''_4, \qquad m = 5, 6, 7 \\
   \MAT{M}_{2*} &\leftarrow (\MAT{M}_{1*} - \MAT{M}_{2*}) / c''_4 \\
   \MAT{M}_{5*} &\leftarrow (\MAT{M}_{2*} - 3 \MAT{M}_{5*}) / c_2 \\
   \MAT{M}_{2*} &\leftarrow \MAT{M}_{2*} + 2 (1 - 4 c_2) \MAT{M}_{5*} \\
   \left[ \begin{array}{c} \MAT{M}_{2*} \\ \MAT{M}_{5*} \end{array} 
\right] &= \left[ \begin{array}{cccc} c_3 (c_3 - c_2) & \, 2 c_2 \, & \, 
c_3 \, & c'_3 c_2 \\ c'_3 & -c_3 & 0 & 0 \end{array} \right]
 \end{align*} These operations do not destroy the proportionality of the 
\FIRST{} and \THIRD{} rows.

The matrix $\MAT{M}$ depends on $c'_4$, $c'_5$, and $c'_6$ in a linear 
way. The \SECOND{} and \NTH{5} rows of the transformed $\MAT{M}$ do 
depend on $c_2$, $c_3$, and $c'_3$ only. The first three columns of 
$\MAT{M}$ form a rank-deficient matrix if
 \begin{align}
   c'_m &= \frac {c'_3 c_3 c_m^2 (c_m - c_2) + 3 c''_m (c_3 - c_2) 
(c_3^2 - 2 c'_3)} {c_3^3 (c_m - c_2) - c_2 (c_m - c_3) (c_3^2 - 2 c'_3)
}, \quad m = 4, 5, 6 \label{cpm}
 \end{align} Note that the expression (\ref{cpm}) for $c'_m$ is valid 
for any $1 \le m \le 7$. Indeed, for $m = 1$ and $m = 2$ the 
eq.~(\ref{cpm}) gives $c'_1 = c'_2 = 0$ due to $c_1 = c''_1 = 0$ and 
$c_m - c_2 = c''_2 = 0$, respectively. For $m = 3$ due to $c_m - c_3 = 
c''_3 = 0$ the eq.~(\ref{cpm}) is reduced to a tautology $c'_3 = c'_3$. 
For $m = 7$, as $c_7 = 1$ and $c''_7 = 1 / 6$, the expression gives 
$c'_7 = 1 / 2$. Also if $c'_3 = c_3^2 / 2$, then the eq.~(\ref{cpm}) 
gives $c'_m = c_m^2 / 2$ whenever $c_m \ne c_2$.

Below $\MAT{M} = \bigl[ m_{ i  j} \bigr]$ stands for 
the already transformed matrix. It is of rank $2$, as $m_{51} m_{24} - 
m_{21} m_{54} = c'^2_3 c_2 \ne 0$. Since $m_{53} = m_{54} = 0$ and 
$m_{23} = c_3$, $m_{24} = c'_3 c_2$, the following linear combinations 
$q_i = c'_3 c_2 m_{ i  3} - c_3 m_{ i 
 4}$, where $i = 1$, $3$, $4$, should be equal to zero. The 
equation $q_1 = 0$ is linear in $c''_4$, with the solution
 \begin{align*}
  c'_4 &= \frac {c'_3 c_4^2 (c_4 - c_2) \bigl( c_3^2 (c_3 - c_2) + c'_3 
(3 c_2 - 2 c_3) \bigr)} {c_4 \bigl( 2 c'^2_3 c_2 + c_3 (c_3 - c_2)^2 
(c_3^2 - 2 c'_3) \bigr) - c'_3 c_2 c_3 \bigl( 2 c'_3 - c_3 (c_3 - c_2) 
\bigr)} \\
  c''_4 &= \frac {c'^2_3 c_2 c_4^2 (c_4 - c_2) (c_4 - c_3)} {c_4 \bigl( 
2 c'^2_3 c_2 + c_3 (c_3 - c_2)^2 (c_3^2 - 2 c'_3) \bigr) - c'_3 c_2 c_3 
\bigl( 2 c'_3 - c_3 (c_3 - c_2) \bigr)}
  \end{align*} The numerator of $q_3$ is bilinear in $c''_4$ and 
$c''_5$. With $c''_4$ being already set, the variable $c''_5$ is 
determined from effectively a linear equation $q_3 = 0$. This results in 
$q_4 = 0$ and $\MAT{M}$ being of rank $2$, also $\textrm{rank}\, \bigl[ 
\, \MAT{M}_{1*}^{\textrm{T}} ~~ \MAT{M}_{3*}^{\textrm{T}} \, \bigr] = 
1$. The variable $c''_6$ is found from a linear equation $\mbox{``the 
numerator of $\bigl( \VEC{b}^{\textrm{T}} (\VEC{c} * \VEC{c} * \VEC{c} * 
\VEC{c}) - 1 / 5 \bigr)$''} = 0$.

The expressions for $c'_5$, $c''_5$, and $c'''_5$ (and especially for 
$c'_6$, $c''_6$, $c'''_6$, $a_{65}$, $b_6$, $d_5$, $d_6$, and $d_7$) are 
too bulky to be included in this paper.

Some combinations of the variables can be written in a relatively 
compact form. For example, here is the expression for stability function 
$R(z)$ that determines the region of absolute stability (see, \EG, 
\citep[sec.~238]{But16}, \citep[sec.~4.4]{AsPe98}):\footnote{Compare 
with \citep[eq.~(3.2)]{DoPr80} and \citep[eq.~(16)]{PaPa96} (the latter 
contains a sign error), where $c'_3 = c_3^2 / 2$.}
 \begin{gather*}
  R(z) = 1 + z  \VEC{b}^{\textrm{T}} (\MAT{I} - z  
\MAT{A})^{-1} \VEC{1} = 1 + z + {\TS\frac1{2}} z^2 + {\TS\frac1{6}} z^3 
+ {\TS\frac1{24}} z^4 + {\TS\frac1{120}} z^5 + b_6 a_{65} c'''_5 z^6 \\
   b_6 a_{65} c'''_5 = - \frac {\gamma_{A c^2} + 5 c''_4 \eta_{A c^2}} 
{120  \GAMMAC{4}{2}} = \frac{c_4}{120} \biggl( 1 - \frac {5 
c'_3 c_2} {2 c'_3 - c_3 (c_3 - c_2)} \biggr)
 \end{gather*} Here $\MAT{I}$ is the identity matrix.

\section{Five families of embedded pairs} \label{five_families}

With $\VEC{c}'$, $\VEC{c}''$, $\VEC{c}'''$, $a_{65}$, and $b_6$ 
expressed through $\VEC{c}$ and $c'_3$, and all but four order 
conditions being met; the three conditions $\VEC{b}^{\textrm{T}} \bigl[ 
\, (\VEC{c}' * \VEC{c} * \VEC{c}) ~~~ \bigl( \bigl( \MAT{A} (\VEC{c} * 
\VEC{c}) \bigr) * \VEC{c} \bigr) ~~~ (\VEC{c}'' * \VEC{c}) \, \bigr] = 
\bigl[ \, \frac{1}{10} ~~ \frac{1}{15} ~~ \frac{1}{30} \, \bigr]$ are 
satisfied when\footnote{They are also satisfied when $(c_3 - c_2) (c_3^2 
- 2 c'_3) = 0$, $c_3^2 (c_3 - c_2) - c'_3 c_2 = 0$, and $c_3^2 (c_3 - 
c_2) + 2 c'_3 c_2 = 0$, respectively. Not satisfying any of 
eq.~(\ref{c4_eq_something}) and eq.~(\ref{c6_eq_1}) would imply $c'_3 
c_2 = 0$ then.}${}^{,}$\footnote{\;Compare with \citep[eq.~(20)]{Feh69}, 
\cite[eq.~(3.3)]{DoPr80}, \citep[p.~1173, Corollary~1]{PaPa96}.}
 \begin{gather}
   c_4 = \frac {c'_3 c_2 \bigl( 2 c'_3 - c_3 (c_3 - c_2) \bigr)} {\bigl( 
2 c'_3 (1 - 2 c_2)- c_3 (c_3 - c_2) \bigr)^2 + 4 c'^2_3 c_2^2} 
\label{c4_eq_something}
 \end{gather} or \begin{gather} c_6 = 1 \label{c6_eq_1} \end{gather} If 
the node $c_4$ is chosen according to eq.~(\ref{c4_eq_something}), then 
the \THIRD{} and \NTH{4} rows of the untransformed matrix $\MAT{M}$ are 
proportional to each other, which results in $d_7 = 0$ and effectively a 
pair of $6$-stage Runge--Kutta methods. The last remaining order 
condition $\VEC{b}^{\textrm{T}} (\VEC{c}' * \VEC{c}') = \frac{1}{20}$ is 
met in three cases:
 \begin{align*}
   \textrm{type~A:}& \quad c'_3 = c_3^2 / 2 \\
   \textrm{type~B:}& \quad c'_3 = 3 (c_3 - c_2) (c_2 + c_3 - 4 c_2 c_3) 
\mathop{\big/}  2 \bigl( 3 - 12 c_2 + 10 c_2^2 \bigr) \\ 
   \textrm{type~C:}& \quad \bigl[ 3 Z(12, 15, 20) - 3 (c_2 + c_3) Z(33, 
40, 50) + 2 c_2 c_3 Z(138, 165, 200) \bigr] {} \\ {} & \quad {} \cdot 
c_3^2 (c_3 - c_2)^2 - \bigl[ (12 + 50 c_2^2) \bigl( Z(12, 15, 20) - c_3 
Z(33, 40, 50) \bigr) {} \\ {} & \quad {} - 3 c_2 Z(207, 260, 
350) + 2 c_2 c_3 Z(852, 1035, 1300) \bigr] c'_3 c_3 (c_3 - c_2) {} \\ {} 
& \quad {} + \bigl[ (2 + 10 c_2 c_3) Z(12, 15, 20) - 15 c_2 
Z(3, 4, 6) - 2 c_3 Z(33, 40, 50) \bigr] {} \\ {} & \quad {} 
\cdot 2 (3 - 12 c_2 + 10 c_2^2)
 c'^2_3 = 0 \end{align*} where $Z(\alpha_0, \alpha_1, \alpha_2) = 
\alpha_0 - \alpha_1 (c_5 + c_6) + \alpha_2 c_5 c_6$. The left-hand side 
in the condition for embedded pairs of type~C is bilinear in $c_5$ and 
$c_6$. Formulas for embedded pairs of type~A are available in 
Appendix~\ref{formulas_A} (see also \citep[app.]{PaPa96}); for pairs of 
type~B see Appendix~\ref{formulas_B}.

If $c_6 = 1$, then the condition $\VEC{b}^{\textrm{T}} (\VEC{c}' * 
\VEC{c}') = \frac{1}{20}$ is met in two cases:
 \begin{align*}
   \textrm{type~A${}'$:}& \quad c'_3 = c_3^2 / 2 \\
   \textrm{type~B${}'$:}& \quad \textrm{a bulky expression (which is a 
polynomial of $c_2$, $c_3$, $c_4$, $c_5$, and $c'_3$} \\ {} & \quad 
\textrm{with degrees $8$, $17$, $3$, $2$, and $8$, respectively) is 
equal to zero}{\protect\footnotemark}
 \end{align*}
 \footnotetext{\;If instead of $c'_3$ the variable $g'_3 = 
c'_3 / c_3 (c_3 - c_2)$ is used, then the degrees are $6$, $2$, $3$, 
$2$, and $8$.}Formulas for embedded pairs of type~A${}'$ are available 
in \citep[app.]{PaPa96}. For pairs of type~B${}'$ the expressions are 
simplified in the cases $c_3 = 0$ (see 
Appendix~\ref{formulas_Bp_c3_eq_0}) and $c_3 = c_2$ (see 
Appendix~\ref{formulas_Bp_c3_eq_c2}).

\begin{figure} \centerline{\includegraphics{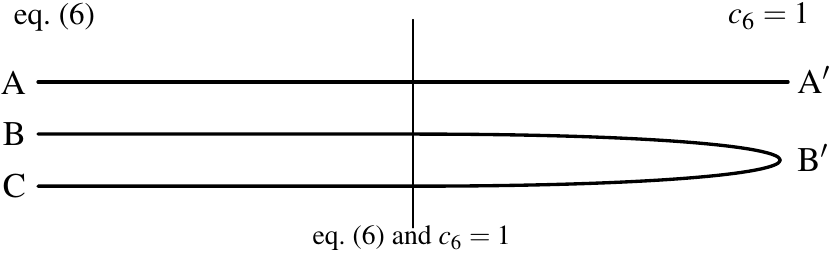} 
\vspace{-2pt}} \caption{Schematic depiction of the five families. The 
left half contains non-FSAL pairs of $6$-stage methods, on the right are 
pairs of $7$-stage methods with FSAL property.} \label{families} 
\end{figure}

\begin{figure} \centerline{\includegraphics{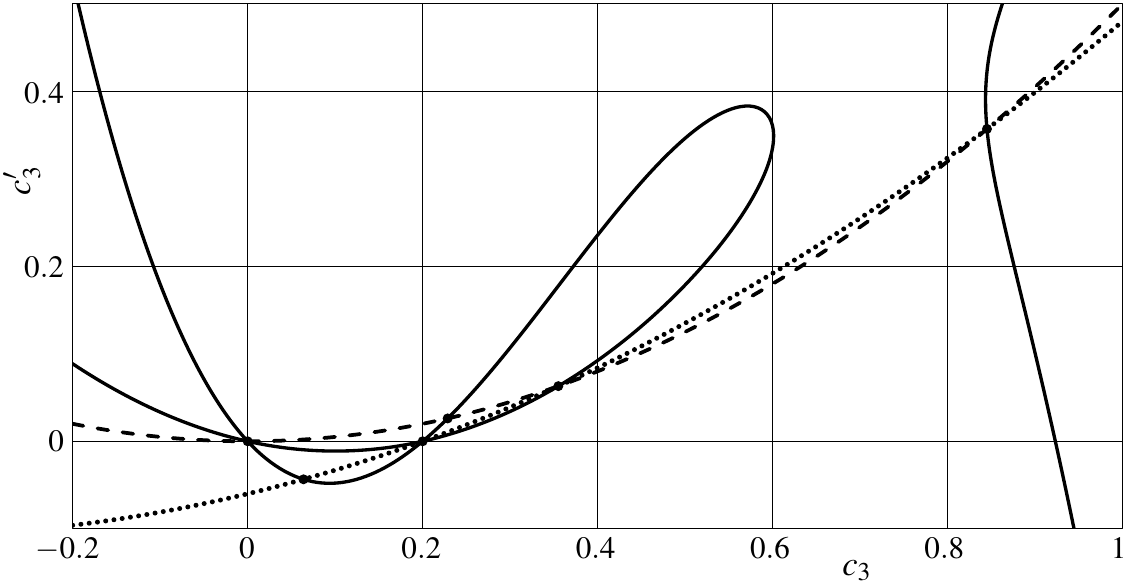} 
\vspace{-4pt}} \caption{A two-dimensional cut through the 
six-dimensional space $(c_2, c_3, c_4, c_5, c_6, 
c'_3)$.{\protect\footnotemark} Here $c_2 = 1 / 5$ and $c_5 = 
4 / 5$. The nodes $c_4$ and $c_6$ are set according to the 
eqs.~(\ref{c4_eq_something}) and (\ref{c6_eq_1}), respectively. The 
dashed, dotted, and solid curves correspond to pairs of type~A, B, and 
C, respectively. The equations for the curves are (A) $c'_3 = c_3^2 / 
2$, (B) $c'_3 = 3 (5 c_3 - 1) (1 + c_3) / 50$, and (C) $c'_3 = c_3 (5 
c_3 - 1) \bigl[ 13 - 12 c_3 \pm (73 - 208 c_3 + 144 c_3^2)^{1 / 2} 
\bigr] / 20 $. All the three curves intersect at $c_3 = \bigl( 6 \pm 
6^{1 / 2} \bigr) / 10$, or when $3 - 12 c_3 + 10 c_3^2 = 0$. The type~C 
curve intersects twice with the ones of type~A and B at $(c_3, c'_3) = 
(0, 0)$ and $(c_3, c'_3) = (c_2, 0)$, respectively. (At these four 
points some of the matrix elements of $\MAT{A}$ are infinite, so they do 
not correspond to any embedded pairs.) The structure of intersections 
stays the same even when only one of the eqs.~(\ref{c4_eq_something}) 
and (\ref{c6_eq_1}) is satisfied.} \label{c3cp3} \end{figure} 
\footnotetext{\;Pairs form a set of codimension $2$ in the 
$6$-dimensional space $(\VEC{c}, c'_3)$. In Figure~\ref{c3cp3} the 
curves in the cut have codimension $1$, as at least one of the 
eqs.~(\ref{c4_eq_something}) and (\ref{c6_eq_1}) (in fact, both) is 
satisfied.}

The connections between pairs of types~A, A${}'$ (that are derived in 
\citep{PaPa96}), B, C, and B${}'$ are shown in Figures~\ref{families} 
and \ref{c3cp3}. The new pairs presented in this paper are listed in the 
lower half of Table~\ref{table_methods}. They were selected by generally 
following the perceptive reasoning in \citep[p.~785]{Ver78}, 
\citep[sec.~3]{DoPr80}, \citep[p.~20]{BoSh96}. As in \citep{DoPr80}, the 
local error was estimated through the $\ell^2$-norms of elementary 
differentials vectors: \begin{align*} T_{p}^2 = 
\sum_{\mathclap{\text{rooted~trees~}\textrm{t}\text{~of~order~}p}} 
\tau^2(\textrm{t}), \qquad \tau(\textrm{t}) = 
\frac{1}{\sigma(\textrm{t})} \biggl( \VEC{b}^{\textrm{T}} 
\VEC{\Phi}(\textrm{t}) - \frac{1}{\textrm{t}!} \biggr) \end{align*} Here 
$\sigma(\textrm{t})$ is the order of the symmetry group of the tree 
$\textrm{t}$ (see, \EG, \citep[p.~154]{But16}, \citep[p.~58]{But21}). 
First, the local error $T_6$ was minimized with inequality constraints 
$\max_{ij} \vert a_{ij} \vert < M$ (for some limit $M$) and $\min_j b_j
> -3$. Then the pair were chosen close to the optimum, with 
representation of coefficients $a_{ij}$ requiring a small number of 
digits. The pair of type~A${}'$ in Table~\ref{method_similar_to_Tsit5} 
was constructed to be a close analogue of \citep[tab.~1]{Tsi11} pair, 
which is of type~B${}'$.

\begin{figure} \centerline{\includegraphics{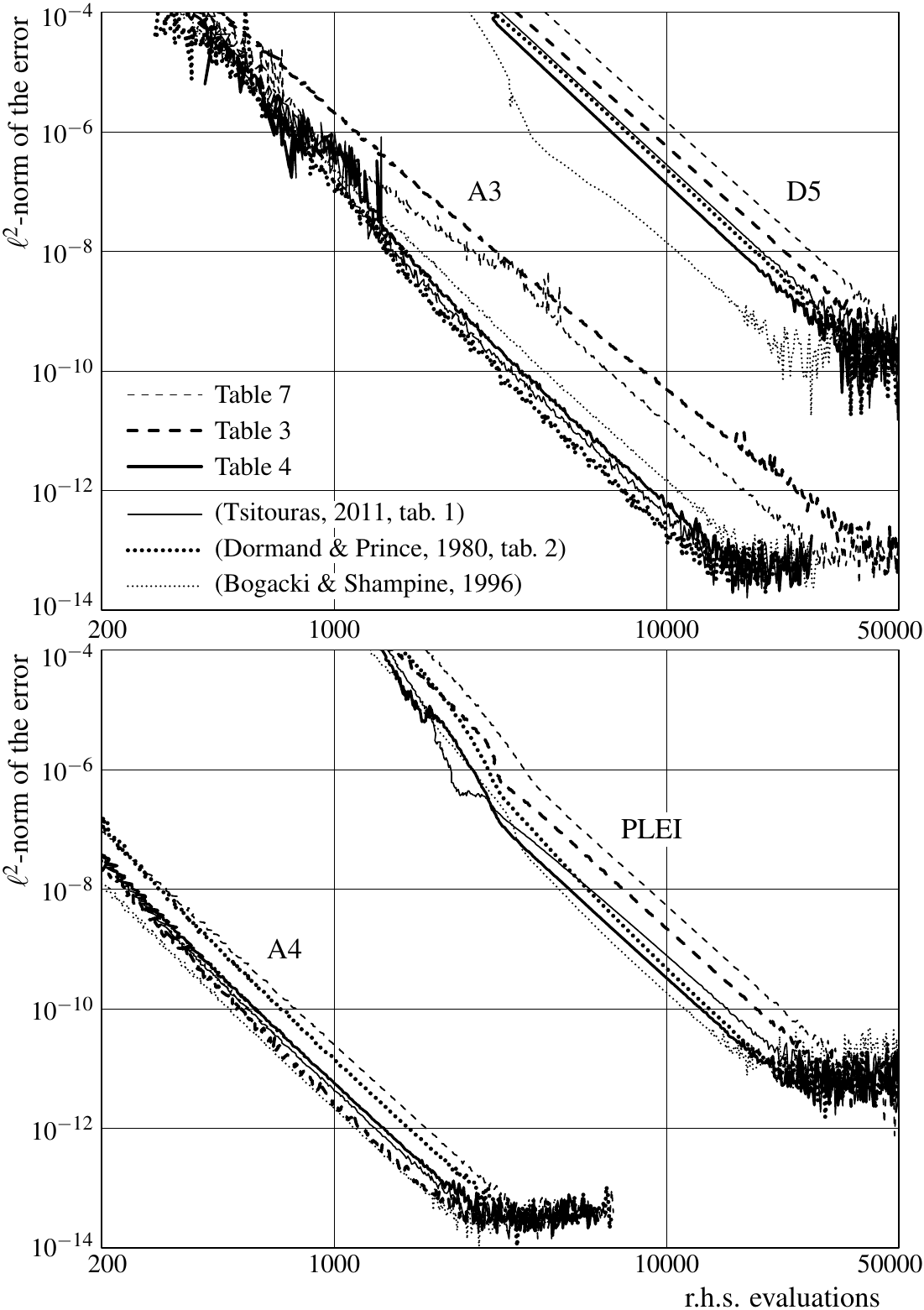} 
\vspace{-4pt}} \caption{Efficiency curves for problems A3, A4 
\citep[p.~617]{DETEST}, D5 \citep[p.~620]{DETEST}, and PLEI 
\citep[p.~245]{HNW93}: the pair in Table~\ref{type_B_pair} (dashed 
curve), Table~\ref{sqrt4054} (thin dashed curve), 
Table~\ref{method_similar_to_Tsit5} (solid curve), and 
\citep[tab.~1]{Tsi11} (thin solid curve), \citep[tab.~2]{DoPr80} (dotted 
curve), and \citep{BoSh96} (thin dotted curve) pairs. The adaptive step 
size scheme $h \leftarrow 0.9 h \bigl( \mbox{ATOL} / E \bigr){}^{1 / 5}$ 
was used. (The starting step size $h_0 = 10^{-6}$ was swiftly corrected 
by the adaptive step size control.) Here $\mbox{ATOL}$ is the absolute 
error tolerance, and $E$ is the $\ell^2$-norm of the difference vector 
between the two solutions within a pair. The steps with $E > 
\mbox{ATOL}$ were rejected, but they were still contributing to the 
number of the r.h.s.~evaluations. For A3, A4, and D5 problems the 
maximal value of the $\ell^2$-norm of the error $\Vert 
\VEC{\tilde{x}}(t) - \VEC{x}(t) \Vert_2$ along the whole trajectory $0 
\le t \le 20$ is plotted. For PLEI the $\ell^2$-norm of the error was 
measured at the end of the integration interval $t = 3$, using only $14$ 
components of $\VEC{x}$ that correspond to the coordinates of the 
stars.} \label{tests} \end{figure}

The efficiency curves or work-precision diagrams of six pairs (three 
from literature and three new ones) are shown in Figure~\ref{tests}. The 
performance of type B${}'$ pairs in Table~\ref{sqrt4054} is the worst. 
With the exception of problem~A4, the type~B pair in 
Table~\ref{type_B_pair} is the second-worst. The efficiency of 
\citep{BoSh96} pair shows the potential benefit of adding a stage. The 
performance of the three other pairs, \citep[tab.~2]{DoPr80}, 
\citep[tab.~1]{Tsi11}, and Table~\ref{method_similar_to_Tsit5}, is 
comparable. (See \citep[tab.~2]{Tsi11} for the comparison of 
\citep[tab.~2]{DoPr80} and \citep[tab.~1]{Tsi11} pairs on all the $25$ 
problems from \citep{DETEST}.)

 \begin{table} \begin{tabular}{r|llllc}
     & $\; 10^4 \times T_6$ & $\; 10^3 \times T_7$ & {$\, 
\mbox{max}_{ij} \vert a_{ij} \vert$} & {$\;\; \mbox{min}_j b_j \STRUT$} 
& $b_6 a_{65} c'''_5$ \\
     \hline
     \citep[tab.~III]{Feh69} & $33.557...$ & $6.7653...$ & 
$\phantom{0}8$ & $-0.18$ & $1 / 2080$ \\
     \citep[\mbox{eq.~(5)}]{CaKa90} & $\phantom{0}9.4828...$ & 
$1.3689...$ & $\phantom{0}2.5925...$ & $\phantom{-}0.0978...$ & $1 / 
800$ \\
     \citep[\mbox{tab.~2}]{DoPr80} & $\phantom{0}3.9908...$ & 
$3.9557...$ & $11.595...$ & $-0.3223...$ & $1 / 600$ \\
     \citep[\mbox{tab.~1}]{Tsi11} & $\phantom{0}1.3851...$ & $2.1124...$ 
& $12.920...$ & $-3.2900...$ & $1 / 698...$ \\
     \citep{BoSh96} & $\phantom{0}0.2216...$ & $0.2126...$ & 
$\phantom{0}1.1637... $ & $\phantom{-}0.0086...$ & N/A \\
     \hline
     type~B,~Table~\ref{type_B_pair} & $\phantom{0}8.9041...$ & 
$1.2159...$ & $\phantom{0}1.6014...$ & $-0.3077...$ & $7 / 5440$ \\
     type~A${}'$,~Table~\ref{method_similar_to_Tsit5} & 
$\phantom{0}1.2239...$ & $1.9225...$ & $10.435...$ & $-2.9044...$ & $3 / 
2080$ \\
     type~B${}'$,~$c_3 = 0$, Table~\ref{method_c3_eq_0} & 
$\phantom{0}7.6950...$ & $1.6029...$ & $\phantom{0}3.1358...$ & 
$-0.0182...$ & $1 / 720$ \\
     type~B${}'$,~$c_3 = c_2$, Table~\ref{method_c3_eq_c2} & $18.132...$ 
& $2.7565...$ & $19.285...$ & $\phantom{-}0.0416...$ & $1 / 960$ \\
     type~B${}'$,~Table~\ref{sqrt4054} & $\phantom{0}5.6328...$ & 
$1.0199...$ & $\phantom{0}5.8955...$ & $-0.1160...$ & $1 / 600$ 
\end{tabular}
  \caption{A comparison of ten embedded $(4, 5)$ pairs. The first five 
are from the literature. The Fehlberg (also available in 
\citep[tab.~1]{Feh70}), Cash--Karp, Dormand--Prince, and Tsitouras pairs 
are of type~A, A, A${}'$, and B${}'$, respectively. The $\min_j b_j$ 
column shows the minimal value of a non-zero weight. The quantity $b_6 
a_{65} c'''_5$ is the coefficient at $z^6$ in the stability function 
$R(z)$. The stability region is most extended when its value is around 
$1 / 1280$ \citep[fig.~2]{Law66}. The Bogacki--Shampine pair is non-FSAL 
and uses $7$ stages, so its absolute stability region is not determined 
by the value of $b_6 a_{65} c'''_5$.} \label{table_methods} \end{table}

 \begin{table} \begin{tabular}{c|cccccc}
     $0\STRUT$ \\
     $\frac{1}{6}\STRUT$  & $\phantom{-}\frac{1}{6}$ \\
     $\frac{7}{32}\STRUT$ & $\phantom{-}\frac{67}{512}$ & 
$\phantom{-}\frac{45}{512}$ \\
     $\frac{33}{68}\STRUT$ & $\phantom{-}\frac{224787}{903992}$ & 
$-\frac{1233765}{903992}$ & $\phantom{-}\frac{180960}{112999}$ \\
     $\frac{3}{4}\STRUT$ & $\phantom{-}\frac{921}{3496}$ & 
$-\frac{552447}{1136200}$ & $\phantom{-}\frac{125664}{316825}$ & 
$\phantom{-}\frac{103173}{179075}$ \\
     $\frac{7}{8}\STRUT$ & $\phantom{-}\frac{13}{13984}$ & 
$-\frac{5604237}{49992800}$ & $\phantom{-}\frac{2246076}{3485075}$ & 
$-\frac{1822723}{189103200}$ & $\phantom{-}\frac{371}{1056}$ \\
     \hline
     $\STRUT$ & $\phantom{-}\frac{1}{9}$ & $-\frac{59508}{193375}$ & 
$\phantom{-}\frac{2281472}{3882375}$ & $\phantom{-}\frac{1920983}{7492875}$ 
& $\phantom{-}\frac{437}{5355}$ & $\phantom{-}\frac{76912}{283815}$ \\
     $\STRUT$ & $\phantom{-}0$ & $\phantom{-}\frac{2349}{700}$ & 
$-\frac{832}{175}$ & $\phantom{-}\frac{83521}{31800}$ & 
$-\frac{377}{168}$ & $\phantom{-}\frac{377}{371}$
 \end{tabular} \caption{An embedded pair of type~B.} \label{type_B_pair} 
\end{table}

 \begin{table} \begin{tabular}{c|ccccccc}
     $0\STRUT$  \\
     $\frac{1}{5}\STRUT$ & $\phantom{-}\frac{1}{5}$ \\
     $\frac{21}{65}\STRUT$ & $\phantom{-}\frac{21}{338}$ & 
$\phantom{-}\frac{441}{1690}$ \\
     $\frac{9}{10}\STRUT$ & $\phantom{-}\frac{639}{392}$ & 
$-\frac{729}{140}$ & $\phantom{-}\frac{1755}{392}$ \\
     $\frac{39}{40}\STRUT$ & $\phantom{-}\frac{4878991}{1693440}$ & 
$-\frac{16601}{1792}$ & $\phantom{-}\frac{210067}{28224}$ & 
$-\frac{1469}{17280}$ \\
     $1\STRUT$ & $\phantom{-}\frac{13759919}{4230954}$ & 
$-\frac{2995}{287}$ & $\phantom{-}\frac{507312091}{61294590}$ & 
$-\frac{22}{405}$ & $-\frac{7040}{180687}$ \\
     $1\STRUT$ & $\phantom{-}\frac{1441}{14742}$ & $\phantom{-}0$ & 
$\phantom{-}\frac{114244}{234927}$ & $\phantom{-}\frac{118}{81}$ & 
$-\frac{12800}{4407}$ & $\phantom{-}\frac{41}{22}$ \\
     \hline
     $\STRUT$ & $\phantom{-}\frac{1441}{14742}$ & $\phantom{-}0$ & 
$\phantom{-}\frac{114244}{234927}$ & $\phantom{-}\frac{118}{81}$ & 
$-\frac{12800}{4407}$ & $\phantom{-}\frac{41}{22}$ \\
     $\STRUT$ & $-\frac{1}{273}$ & $\phantom{-}0$ & 
$\phantom{-}\frac{2197}{174020}$ & $-\frac{4}{15}$ & 
$\phantom{-}\frac{1280}{1469}$ & $-\frac{33743}{52712}$ & 
$\phantom{-}\frac{127}{4792}$ \\
     \hline
     $\theta^{\phantom{2}}\STRUT$ & $\phantom{-}1$ & $\phantom{-}0$ & 
$\phantom{-}0$ & $\phantom{-}0$ & $\phantom{-}0$ & $\phantom{-}0$ & 
$\phantom{-}0$ \\
     $\theta^2\STRUT$ & $-\frac{4489}{1638}$ & $\phantom{-}0$ &
$\phantom{-}\frac{35152}{8701}$ & $-\frac{118}{9}$ & 
$\phantom{-}\frac{48000}{1469}$ & $-\frac{246}{11}$ & 
$\phantom{-}\frac{3}{2}$ \\
     $\theta^3\STRUT$ & $\phantom{-}\frac{21170}{7371}$ & $\phantom{-}0$ 
& $-\frac{1441232}{234927}$ & $\phantom{-}\frac{2596}{81}$ & 
$-\frac{339200}{4407}$ & $\phantom{-}\frac{574}{11}$ & $-4$ \\
     $\theta^4\STRUT$ & $-\frac{2540}{2457}$ & $\phantom{-}0$ & 
$\phantom{-}\frac{202124}{78309}$ & $-\frac{472}{27}$ & 
$\phantom{-}\frac{60800}{1469}$ & $-\frac{615}{22}$ & 
$\phantom{-}\frac{5}{2}$
   \end{tabular} \caption{An embedded pair of type~A${}'$ which is 
structurally similar to the \citep[tab.~1]{Tsi11} pair of type~B${}'$ 
(in the latter one should read $\tilde{b}_7 = -\frac{1}{66}$, also the 
presented vector $\tilde{\VEC{b}}$ is the difference vector $\VEC{d}$). 
The last $4$ rows contain coefficients for the \NTH{4} order 
continuously differential interpolant $\VEC{\tilde{x}}(t + \theta h) = 
\VEC{x}(t) + h \sum_j \beta_j(\theta) \VEC{F}_j = \VEC{x}(t) + h \sum_j 
\VEC{F}_j \sum_k \beta_{kj} \theta^k$, \EG, $\beta_7(\theta) = 
\frac{3}{2} \theta^2 - 4 \theta^3 + \frac{5}{2} \theta^4$.} 
\label{method_similar_to_Tsit5} \end{table}

 \begin{table} \begin{tabular}{c|ccccccc}
     $0\STRUT$ \\
     $\frac{4}{15}\STRUT$ & $\phantom{-}\frac{4}{15}$ \\
     $0\STRUT$ & $\phantom{-}\frac{6}{7}$ & $-\frac{6}{7}$ \\
     $\frac{1}{2}\STRUT$ & $-\frac{11}{384}$ & $\phantom{-}\frac{21}{32}$ & 
$-\frac{49}{384}$ \\
     $\frac{4}{5}\STRUT$ & $\phantom{-}\frac{4}{75}$ & $-\frac{6}{35}$ & 
$\phantom{-}\frac{14}{75}$ & $\phantom{-}\frac{128}{175}$ \\
     $1\STRUT$ & $\phantom{-}\frac{81}{224}$ & $\phantom{-}\frac{4917}{1568}$ 
& $-\frac{33}{32}$ & $-\frac{132}{49}$ & $\phantom{-}\frac{275}{224}$ \\
     $1\STRUT$ & $\phantom{-}\frac{41}{384}$ & $\phantom{-}\frac{3375}{9856}$ 
& $-\frac{7}{384}$ & $\phantom{-}\frac{4}{21}$ & 
$\phantom{-}\frac{125}{384}$ & $\phantom{-}\frac{7}{132}$ \\
     \hline
     $\STRUT$ & $\phantom{-}\frac{41}{384}$ & $\phantom{-}\frac{3375}{9856}$ & 
$-\frac{7}{384}$ & $\phantom{-}\frac{4}{21}$ & 
$\phantom{-}\frac{125}{384}$ & $\phantom{-}\frac{7}{132}$ \\
      $\STRUT$ & $\phantom{-}\frac{1}{40}$ & $\phantom{-}\frac{405}{616}$ & 
$-\frac{7}{40}$ & $-\frac{32}{35}$ & $\phantom{-}\frac{5}{8}$ & 
$-\frac{56}{55}$ & $\phantom{-}\frac{4}{5}$
 \end{tabular} \caption{An embedded pair of type~B${}'$ with $c_3 = 0$.} 
\label{method_c3_eq_0} \end{table}

 \begin{table} \begin{tabular}{c|ccccccc}
     $0\STRUT$ \\
     $\frac{1}{4}\STRUT$ & $\phantom{-}\frac{1}{4}$ \\
     $\frac{1}{4}\STRUT$ & $-\frac{11}{20}$ & $\phantom{-}\frac{4}{5}$ \\
     $\frac{1}{3}\STRUT$ & $\phantom{-}\frac{1}{9}$ & 
$\phantom{-}\frac{43}{216}$ & $\phantom{-}\frac{5}{216}$ \\
     $\frac{4}{5}\STRUT$ & $\phantom{-}\frac{66}{125}$ & $-\frac{593}{250}$ & 
$-\frac{19}{50}$ & $\phantom{-}\frac{378}{125}$ \\
     $1\STRUT$ & $-\frac{7}{2}$ & $\phantom{-}\frac{151}{8}$ & 
$\phantom{-}\frac{25}{8}$ & $-\frac{135}{7}$ & 
$\phantom{-}\frac{25}{14}$ \\
     $1\STRUT$ & $\phantom{-}\frac{5}{48}$ & $\phantom{-}0$ & $\phantom{-}0$ & 
$\phantom{-}\frac{27}{56}$ & $\phantom{-}\frac{125}{336}$ & 
$\phantom{-}\frac{1}{24}$ \\
     \hline
     $\STRUT$ & $\phantom{-}\frac{5}{48}$ & $\phantom{-}0$ & $\phantom{-}0$ & 
$\phantom{-}\frac{27}{56}$ & $\phantom{-}\frac{125}{336}$ & 
$\phantom{-}\frac{1}{24}$ \\
     $\STRUT$ & $\phantom{-}\frac{11}{8}$ & $\phantom{-}\frac{8}{3}$ & 
$-\frac{40}{3}$ & $\phantom{-}\frac{297}{28}$ & $-\frac{125}{56}$ & 
$-\frac{1}{12}$ & $\phantom{-}1$
 \end{tabular} \caption{An embedded pair of type~B${}'$ with $c_3 = 
c_2$. Although $1 / 5 = c'_3 \ne c_3^2 / 2 = 1 / 32$, the weight $b_3 = 
0$ and $c'_m = c_m^2 / 2$ for $m > 3$. Thus, the Dominant Stage-Order 
(DSO) \citep[eq.~(5)]{Ver14} of the \NTH{5}~order method is equal to $2$. 
As $d_2 \ne 0$, the \NTH{4} order method has $\mbox{DSO} = 1$.} 
\label{method_c3_eq_c2} \end{table}

 \begin{table} \begin{gather*}
   \begin{array}{c|ccc}
     0\STRUT  \\
     \frac{1}{5}\STRUT  & \frac{1}{5} \\
     \frac{1}{4}\STRUT  & \frac{1}{8} & \frac{1}{8} \\
     \frac{3}{5}\STRUT  & \frac{141}{575} & -\frac{228}{115} & 
\frac{1344}{575} \\
     c_5\STRUTT & -\frac{c_5 (860 c_5^3 - 1077 c_5^2 + 379 c_5 - 
48)}{3 (39 c_5 - 5)} & \frac{c_5 (5 c_5 - 1) (1340 c_5^2 - 1367 
c_5 + 277)}{2 (39 c_5 - 5)} & -\frac{16 c_5 (5 c_5 - 1) (4 c_5 - 
1) (73 c_5 - 55)}{7 (39 c_5 - 5)} \\
     1\STRUTT & \frac{113 c_5^2 - 35 c_5 - 40}{c_5 (285 - 319 c_5)} & 
-\frac{4 (2845 c_5^2 - 2999 c_5 + 654)}{(5 c_5 - 1) (285 - 319 c_5)} & 
\frac{384 (168 c_5^2 - 193 c_5 + 52)}{7 (4 c_5 - 1) (285 - 319 c_5)} \\
     1\STRUTT & \frac{31 c_5 - 5}{288 c_5} & -\frac{125 (3 - c_5)}{768 
(5 c_5 - 1)} & \frac{8 (7 c_5 + 3)}{63 (4 c_5 - 1)} \\
     \hline
     \STRUTT & \frac{31 c_5 - 5}{288 c_5} & -\frac{125 (3 - c_5)}{768 (5 
c_5 - 1)} & \frac{8 (7 c_5 + 3)}{63 (4 c_5 - 1)} \\
     \STRUTT & \frac{5 (43 c_5 - 33)}{216 c_5} & \frac{175 (43 c_5 - 
33)}{288 (5 c_5 - 1)} & -\frac{152 (43 c_5 - 33)}{189 (4 c_5 - 1)}
   \end{array} \\
   \begin{array}{c|cccc}
     c_5 &\STRUTT \frac{115 c_5 (5 c_5 - 1) (4 c_5 - 1) (5 c_5 - 
3)}{42 (39 c_5 - 5)} \\
     1\STRUTT & -\frac{460 (35 c_5^2 - 55 c_5 + 22)}{7 (5 c_5 - 3) (285 
- 319 c_5)} & \frac{24 (1 - c_5) (39 c_5 - 5)}{c_5 (5 c_5 - 1) (4 c_5 
- 1) (5 c_5 - 3) (285 - 319 c_5)} \\
     1\STRUTT & \frac{2875 (7 c_5 - 5)}{8064 (5 c_5 - 3)} & \frac{39 
c_5 - 5}{96 c_5 (5 c_5 - 1) (4 c_5 - 1) (5 c_5 - 3) (1 - c_5)} & 
\frac{285 - 319 c_5}{2304 (1 - c_5)} \\
     \hline
     \STRUTT & \frac{2875 (7 c_5 - 5)}{8064 (5 c_5 - 3)} & \frac{39 c_5 
- 5}{96 c_5 (5 c_5 - 1) (4 c_5 - 1) (5 c_5 - 3) (1 - c_5)} & \frac{285 - 
319 c_5}{2304 (1 - c_5)} \\
     \STRUTT & \frac{575 (43 c_5 - 33)}{1512 (5 c_5 - 3)} & 
-\frac{(39 c_5 - 5) (43 c_5 - 33)}{72 c_5 (5 c_5 - 1) (4 c_5 - 1) (5 c_5 
- 3) (1 - c_5)} & -\frac{5 (285 - 319 c_5)}{864 (1 - c_5)} & ~~~~1
   \end{array} \\[6pt]
   {\small \begin{array}{c|rrrrrrr}
     0.000\, \\
     0.200\, & 0.200 \\
     0.250\, & 0.125 & 0.125 \\
     0.600\, & 0.245 & -1.983 & 2.337 \\
     0.814\, & -0.107 & 2.416 & -2.110 & 0.615 \\
     1.000\, & 0.304 & -4.967 & 5.896 & -1.014 & 0.782 \\
     1.000\, & 0.086 & -0.116 & 0.490 & 0.232 & 0.249 & 0.059 \\
     \hline
     & 0.086 & -0.116 & 0.490 & 0.232 & 0.249 & 0.059 \\
     & 0.056 & 0.392 & -0.707 & 0.706 & -0.657 & -0.791 & 
\phantom{-}1.000
   \end{array} }
 \\[-24pt] \phantom{.} \end{gather*} \caption{An embedded pair of 
type~B${}'$. The parameters are $c_2 = 1 / 5$, $c_3 = 1 / 4$, $c'_3 = 1 
/ 40$, and $c_4 = 3 / 5$. All the conditions up to the \NTH{5} order are 
satisfied but $\VEC{b}^{\textrm{T}} (\VEC{c}' * \VEC{c}') - 1 / 20 = 
(289 c_5^2 + 2586 c_5 - 2295) / 6900 (39 c_5 - 5) (285 - 319 c_5) = 0$, 
which leads to $c_5 = 3 \bigl( 8 \sqrt{4054} - 431 \bigr) / 289 = 
0.81351...$ (The other choice $c_5 = -9.76...$ would result in $T_6 = 
0.045...$, $T_7 = 0.30...$, and $\vert a_{52} \vert > 8.9 \times 10^4$.) 
At the bottom is the Butcher tableau rounded to the nearest thousandth.} 
\label{sqrt4054} \end{table}

\section{Conclusions}

In pairs of $7$-stage explicit Runge--Kutta methods, the FSAL property 
implies $c_6 = 1$ and the condition $D(1)$: $\sum_{i} b_i a_{ij} = b_j 
(1 - c_j)$ (see, \EG, \citep[p.~189]{But16}, \citep[eq.~(5.6)]{HNW93}, 
\citep[pp.~173 and 193]{But21}), regardless of whether the simplifying 
assumption is satisfied (type~A${}'$) or not (type~B${}'$). There are 
pairs of $8$-stage methods with FSAL property and $c_7 \ne 1$, \EG, the 
\citep[fig.~3]{OwZe92} pair has $\VEC{c}^{\textrm{T}} = \bigl[ \, 0 ~~ 
\frac{1}{6} ~~ \frac{1}{4} ~~ \frac{1}{2}~~ \frac{1}{2} ~~ \frac{9}{14} 
~~ \frac{7}{8} ~~ 1 \, \bigr]$.

The simplifying assumption $c'_i = \sum_{j} \AIJ \mkern1mu c_{\mkern-1mu 
j} = c_i^2 / 2$, where $i \ne 2$, introduces additional redundancy in 
the order conditions, and the number of free parameters in the families 
of pairs of types A and A${}'$ (that do satisfy the assumption) is the 
same as for B, B${}'$, and C (that do not satisfy the assumption). Not 
assuming the simplifying assumption does not increase the dimension of 
the set of pairs satisfying the order conditions. The pairs of 
types~A${}'$ and B${}'$ form different $4$-dimensional submanifolds of 
the space of matrices $\MAT{A}$, with a $3$-dimensional intersection.

From numerical experiments, the part of type~B${}'$ pairs set that 
contains efficient pairs is close to the set of type~A${}'$ pairs. For 
example, in \citep[tab.~1]{Tsi11} pair the weight $b_2 = \frac{1}{100}$ 
is small, and $a_{32} c_2 / (c_3^2 / 2) = 1.0102...$ is close to $1$. It 
is hard to expect a good pair of type~B${}'$ without a counterpart of 
type~A${}'$.

\medskip

The author is grateful to anonymous reviewers for helpful comments and 
suggestions.

\appendix

\begingroup\allowdisplaybreaks

\section{Formulas for pairs of type~A} \label{formulas_A}

 \begin{align*} \\[-27pt]
    c_4 &= c_3 \mathop{\big/} 2 (1 - 4 c_3 + 5 c_3^2) \\
    c'_m &= c_m^2 / 2 , \qquad m = 3, 4, 5, 6 \\
    c''_m &= c_m (c_m - c_3) (c_3 + c_m - 4 c_3 c_m) \mathop{\big/} 2 
\bigl( 3 - 12 c_3 + 10 c_3^2 \bigr) , \qquad m = 4, 5, 6 \\
    c'''_5 &= c_3 c_5 (c_5 - c_3) (c_5 - c_4) \mathop{\big/} 4 (3 - 12 
c_3 + 10 c_3^2) \\
    g &= 8 c_3 - 15 c_3^2 - 4 c_5 (1 - 4 c_3 + 5 c_3^2) + 2 
c_6 (2 - 13 c_3 + 20 c_3^2) \\
    c'''_6 &= \frac {g \, c_6 (c_6 - c_3) (c_6 - c_4)} {4 (3 - 12 
c_3 + 10 c_3^2) (8 - 15 c_3 - 10 c_5 + 20 c_3 c_5)} \\
    b_6 a_{65} c'''_5 &= c_4 (2 - 5 c_3) / 240 \\
    b_6 c'''_6 &= g / 480 (c_6 - c_5) \bigl( 1 - 4 c_3 + 5 c_3^2 \bigr) 
\\
    b_2 &= d_2 = d_7 = 0 \\
    d_5 c_5 (c_5 &- c_3) (c_5 - c_4) + d_6 c_6 (c_6 - c_3) (c_6 - c_4) = 
0
 \end{align*} Note that $\VEC{c}'$, $\VEC{c}''$, $\VEC{c}'''$, and $b_6$ 
do not depend on $c_2$. As $b_2 = d_2 = 0$, the whole vectors $\VEC{b}$ 
and $\VEC{d}$ do not depend on $c_2$. The coefficients $a_{ij}$ and the 
weights $b_j$, $d_j$ are obtained using formulas in the beginning of 
Section~\ref{rewrite}, \EG, $b_5 = (1 / 120 - b_6 c'''_6) / c'''_5$.

\section{Formulas for pairs of type~B} \label{formulas_B}

 \begin{align*} \\[-27pt]
   g &= (3 - 12 c_2 + 10 c_2^2) (3 - 12 c_3 + 10 c_3^2) + 
15 (c_2 + c_3 - 4 c_2 c_3)^2 \\
   c_4 &= 3 (3 - 10 c_2 c_3) (c_2 + c_3 - 4 c_2 c_3) \mathop{\big/} 2 
g \\
   c'_m &= 3 (c_m - c_2) (c_2 + c_m - 4 c_2 c_m) \mathop{\big/} 2 \bigl( 
3 - 12 c_2 + 10 c_2^2 \bigr) , \qquad m = 3, 4, 5, 6 \\
   h_m &= 3 c_2 + 3 c_3 + 3 c_m - 12 c_2 c_3 - 12 c_2 c_m 
- 12 c_3 c_m + 38 c_2 c_3 c_m \\
   c''_m &= \frac {(c_m - c_2) (c_m - c_3) \, h_m} {2 \bigl( 3 - 12 c_2 
+ 10 c_2^2 \bigr) \bigl( 3 - 12 c_3 + 10 c_3^2 \bigr)} , \qquad m = 4, 
5, 6 \\
   c'''_5 &= \frac {3 (c_5 - c_2) (c_5 - c_3) (c_5 - c_4) (c_2 + c_3 - 4 
c_2 c_3)} {4 \bigl( 3 - 12 c_2 + 10 c_2^2 \bigr) \bigl( 3 - 12 c_3 + 10 
c_3^2 \bigr)} \\
   p &= 24 - 45 c_2 - 45 c_3 + 100 c_2 c_3 - 10 \bigl[ 3 
- 6 c_2 - 6 c_3 + 14 c_2 c_3 \bigr] c_5 \\
   q &= 3 (c_2 + c_3 - 4 c_2 c_3) (24 - 45 c_2 - 45 c_3 + 
100 c_2 c_3) {} \\
   &- \bigl[ 4 \bigl( 3 - 12 c_2 + 10 c_2^2 \bigr) \bigl( 3 - 12 c_3 + 
10 c_3^2 \bigr) + 60 (c_2 + c_3 - 4 c_2 c_3)^2 \bigr] c_5 {} \\
   &+ \bigl[ 4 \bigl( 3 - 12 c_2 + 10 c_2^2 \bigr) \bigl( 3 - 12 c_3 + 
10 c_3^2 \bigr) {} \\
   &{}\qquad\quad {} - 30 (c_2 + c_3 - 4 c_2 c_3) (3 - 8 c_2 - 8 c_3 + 
22 c_2 c_3) \bigr] c_6 \\
   c'''_6 &= \frac {(c_6 - c_2) (c_6 - c_3) (c_6 - c_4) q} 
{4 \bigl( 3 - 12 c_2 + 10 c_2^2 \bigr) \bigl( 3 - 12 c_3 + 10 c_3^2 
\bigr)  p} \\
   b_6 a_{65} c'''_5 &= (c_2 + c_3 - 4 c_2 c_3) (6 - 15 c_2 - 15 c_3 + 
40 c_2 c_3) \mathop{\big/} 160  g \\
   b_6 c'''_6 &= q \mathop{\big/} 480 (c_6 - c_5)  g \\
   b_1 &= 1 / 9 \\
   d_1 &= d_7 = 0 \\
   d_5 (c_5 &- c_2) (c_5 - c_3) (c_5 - c_4) + d_6 (c_6 - c_2) (c_6 - 
c_3) (c_6 - c_4) = 0
 \end{align*}

\section{Formulas for pairs of type~B${}'$, $c_3 = 0$} 
\label{formulas_Bp_c3_eq_0}

 \begin{align*} \\[-27pt]
      c_3 &= 0 \\
      c_6 &= 1 \\
        & \hskip-12pt \begin{tabular}{|c|cc|cc|cc|cc|}
         \hline
 \multirow{2}{*}{$\alpha_{l  m  n}$} & 
\multicolumn{2}{c|}{\scriptsize $l = 0$} & 
\multicolumn{2}{c|}{\scriptsize $l = 1$} & 
\multicolumn{2}{c|}{\scriptsize $l = 2$} & 
\multicolumn{2}{c|}{\scriptsize $l = 3$} \\
         \cline{2-9}
         & {\scriptsize $m = 0$} & {\scriptsize $m = 1$} & 
{\scriptsize $m = 0$} & {\scriptsize $m = 1$} & 
{\scriptsize $m = 0$} & {\scriptsize $m = 1$} & 
{\scriptsize $m = 0$} & {\scriptsize $m = 1$} \\
         \hline
 {\scriptsize $n = 0$} & 144 &  180 & 180 & 228 &  72 &  93 &  9 & 12$$ \\
 {\scriptsize $n = 1$} & 360 &  940 & 512 & 940 & 222 & 366 & 30 & 48 \\
 {\scriptsize $n = 2$} & 200 & 1100 & 340 & 960 & 162 & 360 & 24 & 48$$ \\
         \hline \end{tabular} \\
      g &= 5 (c_2^2 + 4 c_4^2) c_5 (3 - 5 c_5) - c_2 c_4  
\sum_{l = 0}^3 \sum_{m = 0}^1 \sum_{n = 0}^2 (-1)^{ l + m + n} 
\alpha_{l  m  n}  (5 c_2)^l c_4^m c_5^n \\
      p &= 3 - 5 c_2 - 5 c_4 + 10 c_2 c_4 \\
      q &= 12 - 15 c_2 - 15 c_4 - 15 c_5 + 20 c_2 c_4 + 20 c_2 c_5 + 20 
c_4 c_5 - 30 c_2 c_4 c_5 \\
   c'_3 &= 3 g \mathop{\big/} 2  (6 - 15 c_2 - 10 c_5 + 30 c_2 c_5) 
 p  q \\
   c'_4 &= 3 c_4 (c_4 - c_2)  / 2 \\
   c'_5 &= 3 (c_5 - c_2) \bigl( c_5 + c_4 (2 - 5 c_2 - 5 c_5 + 10 c_2 
c_5) \bigr) \mathop{\big/} 2 p \\
   c'_6 &= 3 (1 - c_2) \bigl( 4 - 7 c_4 - 5 c_5 + 5 c_2 c_4 + 10 (1 - 
c_2) c_4 c_5 \bigr) \mathop{\big/} 2  q \\
   c''_m &= c'_m c_m / 3 , \qquad m = 4, 5, 6 \\
   c'''_5 &= c_4 c_5 (c_5 - c_2) (c_5 - c_4) (2 - 5 c_2) / 4 p \\
   c'''_6 &= (1 - c_2) (1 - c_4) (2 - 2 c_4 - 2 c_5 + 5 c_2 c_4) / 4 
 q \\
   b_6 a_{65} c'''_5 &= c_4 (2 - 5 c_2) / 240 \\
   b_6 c'''_6 &= (2 - 2 c_4 - 2 c_5 + 5 c_2 c_4) / 240 (1 - c_5) \\
     d_5 &= p  \bigl( c_2 c_4 + (c_2 - 2 c_4) (3 - 5 c_5) + 
15 c_2 (1 - c_2) c_4 (1 - 2 c_5) \bigr) \\
     d_6 &= q  c_5 (c_5 - c_2) (c_5 - c_4) \bigl( 4 c_4 - 2 
c_2 - 14 c_2 c_4 + 15 c_2^2 c_4 \bigr) \mathop{\big/}  (1 - 
c_2) (1 - c_4) \\
     d_7 &= 15 c_5 (c_5 - c_2) (c_5 - c_4) (1 - c_5) \bigl(c_2 - 2 c_4 + 
8 c_2 c_4 - 10 c_2^2 c_4 \bigr)
 \end{align*}

\section{Formulas for pairs of type~B${}'$, $c_3 = c_2$} 
\label{formulas_Bp_c3_eq_c2}

 \begin{align*} \\[-27pt]
    c_3 &= c_2 \\
    c_4 &= (3 - 5 c_5) \mathop{\big/} 5 (1 - 2 c_5) \\
    c_6 &= 1 \\
    c'_m &= c_m^2 / 2, \qquad m = 4, 5, 6 \\
    c''_4 &= c_4^2 (c_4 - c_2) / 2 \\
    c''_5 &= c_5 (c_5 - c_2) \bigl( c_5 + c_4 (2 - 5 c_2 - 5 c_5 + 10 c_2 
c_5) \bigr) \mathop{\big/} 2 p \\
    c''_6 &= (1 - c_2) \bigl( 4 - 7 c_4 - 5 c_5 + 5 c_2 c_4 + 10 (1 - 
c_2) c_4 c_5 \bigr) \mathop{\big/} 2  q \\
    b_1 &= (1 - 8 c_5 + 10 c_5^2)  / 12 c_5 (5 c_5 - 3) \\
    b_2 &= b_3 = 0 \\
    b_4 &= 125 (2 c_5 - 1)^4 / 12 (5 c_5 - 2) (5 c_5 - 3) (3 - 10 c_5 + 
10 c_5^2) \\
    b_5 &= 1 / 12 c_5 (1 - c_5) (3 - 10 c_5 + 10 c_5^2) \\
    b_6 &= -(3 - 12 c_5 + 10 c_5^2)  / 12 (1 - c_5) (5 c_5 - 
2) \\
    d_5 &= -(1 - c_2) (5 c_5 - 3) (6 - 15 c_2 - 10 c_5 + 30 c_2 c_5) 
\mathop{\big/} 3 c_5 (3 - 10 c_5 + 10 c_5^2) \\
    d_6 &= \bigl( 12 - 52 c_2 + 45 c_2^2 - 5 c_5 (4 - 18 c_2 + 15 c_2^2) 
\bigr) (3 - 12 c_5 + 10 c_5^2) \mathop{\big/} 3 (5 c_5 - 2) \\
    d_7 &= (1 - c_5) \bigl( 6 - 29 c_2 + 30 c_2^2 - 10 c_5 (1 - 5 c_2 + 
5 c_2^2) \bigr) \end{align*} See Appendix~\ref{formulas_Bp_c3_eq_0} for 
the expressions for $p$, $q$, $c'''_5$, $c'''_6$, and $b_6 a_{65} 
c'''_5$. The whole vector $\VEC{b}$ depends on $c_5$ only.

\endgroup

\end{document}